\documentclass[12pt]{article}

\usepackage{amssymb}

\usepackage{amsmath}

\usepackage[mathcal]{euscript}

\usepackage{amsthm}

\newcommand{\XX}{\mathbb{X}}

\newcommand{\e}{\varepsilon}
\newcommand{\cd }{c\`{a}dl\`{a}g }

\newcommand{\EE}{\mathsf{E}}
\newcommand{\PP}{\mathsf{P}}

\newcommand{\JJ}{\mathsf{J}}

\begin{document}

\begin{center}
{\bf \large Necessary and Sufficient Conditions \\ for Convergence  
of First-Rare-Event Times for \vspace{1mm} \\ Perturbed Semi-Markov Processes}
\end{center} 
\vspace{2mm}

\begin{center}
{\large Dmitrii Silvestrov\footnote{Department of Mathematics, Stockholm University, SE-106 81 Stockholm, Sweden. \\ 
Email address: silvestrov@math.su.se}}
\end{center}
\vspace{2mm}

Abstract:
Necessary and sufficient conditions for convergence in distribution of
first-rare-event times and convergence in Skorokhod J-topology of first-rare-event-time processes 
for  perturbed  semi-Markov processes with finite phase space are obtained.  \\



{\bf 1. Introduction}  \\

Random functionals similar with first-rare-event
times are known under different names such as first hitting times, first passage times, absorption times, in theoretical studies,  
and  as  lifetimes, first failure times, extinction times, etc., 
in applications. Limit theorems for such functionals for Markov type processes  have been studied by many researchers.

The case of Markov chains and semi-Markov processes with finite phase
spaces is the most deeply investigated. We refer here to the works
by Simon and Ando (1961),  Kingman (1963), Darroch and Seneta (1965,
1967), Keilson (1966, 1979), Korolyuk (1969), Korolyuk and Turbin
(1970, 1976), Silvestrov (1970, 1971, 1974, 1980, 2014), Anisimov (1971a, 1971b, 1988, 2008), Turbin (1971), Masol and Silvestrov (1972), Zakusilo
(1972a, 1972b),  Kovalenko (1973), Latouch and Louchard (1978),
Shurenkov (1980a, 1980b), Gut and Holst (1984), Brown and Shao
(1987),  Alimov and Shurenkov (1990a, 1990b), Hasin and Haviv
(1992), Asmussen (1994, 2003), Ele\u \i ko and Shurenkov (1995),  Kalashnikov (1997), Kijima (1997), Stewart (1998, 2001),  
Gyllenberg and Silvestrov (1994, 1999, 2000, 2008),  Silvestrov and Drozdenko (2005, 2006a, 2006b), Asmussen and  Albrecher (2010),  Yin and Zhang (2005, 2013), Drozdenko (2007a, 2007b, 2009), Benois, Landim and  Mourragui (2013).

The case of Markov chains and semi-Markov processes with countable
and an arbitrary phase space was treated in works by Gusak and
Korolyuk (1971), Silvestrov (1974, 1980, 1981, 1995, 2000), Korolyuk
and Turbin (1978), Kaplan (1979, 1980), Kovalenko and  Kuznetsov (1981), Aldous (1982), Korolyuk~D.
and Silvestrov (1983, 1984), Kartashov (1987, 1991, 1996, 2013), Anisimov
(1988, 2008), Silvestrov and Velikii (1988), Silvestrov and Abadov (1991,
1993), Motsa and Silvestrov (1996), Korolyuk and Swishchuk (1992), Korolyuk~V.V. and Korolyuk~V.S.
(1999), Koroliuk and  Limnios (2005), Kupsa and Lacroix (2005),   Glynn (2011),  and Serlet (2013). 

We also refer to the books by  Silvestrov (2004) and Gyllenberg
and Silvestrov (2008) and papers by Kovalenko (1994) and Silvestrov D. and Silvestrov S. (2015),  
where one can find comprehensive  bibliographies of works in the area. 

The main features for the most previous results is that they give
sufficient conditions of convergence for such functionals. As a
rule, those conditions involve assumptions, which imply convergence
in distribution for sums of i.i.d random variables distributed as
sojourn times for the semi-Markov process (for every state) to some
infinitely divisible laws plus some ergodicity condition for the
imbedded Markov chain plus condition of vanishing probabilities of
occurring a rare event during one transition step for the semi-Markov
process.

In the context of  necessary and sufficient conditions of convergence  in distribution for first-rare-event-time type functionals, we  would like  to point out the paper by Kovalenko (1965) and  the books by Gnedenko and Korolev (1996)  and Bening and Korolev (2002), where one can find some related results for geometric sums of random variables, and the papers by  Korolyuk and Silvestrov (1983) and Silvestrov and Velikii (1988), where one can find  some related results for first-rare-event-time type functionals defined on Markov chains with arbitrary phase space. 

The results of the present paper relate to the model of perturbed semi-Markov processes with a
finite phase space. Instead of conditions based on ``individual''  distributions of sojourn times,
we use more general and weaker conditions imposed on distributions
sojourn times averaged by stationary distributions of the corresponding 
imbedded Markov chains. Moreover, we show that these conditions are
not only sufficient but also necessary conditions for convergence in
distribution of first-rare-event times and convergence in Skorokhod J-topology of first-rare-event-time processes.
These results give some kind of a ``final solution'' for
limit theorems for first-rare-event times  and first-rare-event-time processes for perturbed semi-Markov process
with a finite phase space. 

The paper generalize and  improve results concerned necessary and sufficient conditions of weak convergence 
for first-rare-event times for semi-Markov process obtained in papers by Silvestrov and Drozdenko (2005, 2006a, 2006b) and 
Drozdenko (2007a, 2007b, 2009). 

First, weaken model  ergodic conditions are imposed on the corresponding 
embedded Markov chains. Second, the above results about weak convergence 
for first-rare-event times are extended, in Theorem 1, to the form of corresponding functional limit theorems for first-rare-event-time processes, 
with necessary and sufficient conditions of convergence. Third, new proofs, based on  general limit theorems for randomly 
stopped stochastic processes, developed and  extensively presented in Silvestrov (2004), are given,  instead of more traditional proofs based on cyclic representations of first-rare-event times if the form of geometrical type random sums. This actually made it possible to get more advanced results in the  form of functional limit theorems. Fourth,  necessary and sufficient conditions of convergence for step-sum reward processes defined on Markov chains are also obtained in the paper. In the context of the present paper,  these results,  formulated in Theorem 2, play an intermediate role. At the same time, they  have their own theoretical and applied values. Finally, we  would like to mention results formulated in Lemmas 1 - 9, which also give some useful supplementary information about asymptotic properties of first-rare-event times and step-sum reward processes. 

We would like to conclude the introduction with the remark that the present paper is a slightly improved version of the 
research report by Silvestrov (2016). \\

{\bf 2. First-rare-event times for perturbed semi-Markov processes} \\

Let $(\eta_{\e, n}, \kappa_{\e, n}, \zeta_{\e, n}), \ n = 0, 1, \ldots$ be, for every $\e \in (0, \e_0]$,  a Markov
renewal process, i.e., a homogenous Markov chain with a phase space $\mathbb{Z} =
 \{1, 2, \ldots, m\} \times [0, \infty) \times \{ 0, 1 \}$, an initial distribution 
$\bar{q}_\e = \langle q_{\e, i} = \PP \{\eta_{\e, 0} = i, \kappa_{\e, 0} = 0, \zeta_{\e, 0} = 0 \} = \PP \{\eta_{\e, 0} = i \}, i \in {\mathbb X} \rangle$  and transition probabilities,
\begin{equation}\label{sad}
\begin{aligned}
& \PP \{ \eta_{\e, n+1} = j,  \kappa_{\e, n+1} \leq t, \zeta_{\e, n+1} = \jmath
/ \eta_{\e, n} = i, \xi_{\e, n} = s, \zeta_{\e, n} = \imath \}   \\
& \quad = \PP \{ \eta_{\e, n+1} = j,  \kappa_{\e, n+1} \leq t, \zeta_{\e, n+1}
= \jmath / \eta_{\e, n} = i \}  \\
& \quad = Q_{\e, ij}(t, \jmath), \ i, j \in \mathbb{X}, \ s, t \geq 0, \ \imath, \jmath = 0, 1.
\end{aligned}
\end{equation}

As is known, the first component $\eta_{\e, n}$ of the above  Markov renewal
process is also a homogenous Markov chain, with the phase space $\mathbb{X} = \{1, 2, \ldots, m\}$, the initial distribution $\bar{q}_\e = \langle q_{\e, i} = 
\PP \{\eta_{\e, 0} = i \}, i \in {\mathbb X} \rangle$  and the transition probabilities,
\begin{equation}\label{proba}
p_{\e, ij} = Q_{\e, ij}(+ \infty, 0) + Q_{\e, ij}(+ \infty, 1),  \, i , j \in \mathbb{X}.
\end{equation}

Also, the  random sequence $(\eta_{\e, n}, \zeta_{\e, n}), \ n = 0, 1, \ldots$ is a Markov renewal process with the phase space $\mathbb{X} \times \{0, 1 \}$, the initial distribution 
$\bar{q}_\e = \langle q_{\e, i} = \PP \{\eta_{\e, 0} = i, \zeta_{\e, 0} = 0 \} = \PP \{\eta_{\e, 0} = i \}, i \in {\mathbb X} \rangle$  and the transition probabilities,
\begin{equation}\label{probak}
p_{\e, ij, \jmath} = Q_{\e, ij}(+ \infty, \jmath),  \, i , j \in \mathbb{X}, \jmath = 0, 1.
\end{equation}

Random variables $\kappa_{\e, n}, n = 1, 2, \ldots$ can be interpreted as sojourn times and random variables 
$\tau_{\e, n} =  \kappa_{\e, 1}  + \cdots +  \kappa_{\e, n}, n = 1, 2, \ldots, \tau_{\e, 0} = 0$ as moments of jumps  for a 
semi-Markov process  $\eta_\e(t), t \geq 0$  defined by the following relation, 
\begin{equation}\label{semi}
\eta_\e(t) = \eta_{\e, n} \quad \mbox{for} \quad \tau_{\e,n} \leq t
<\tau_{\e, n+1}, \  n = 0, 1, \ldots,
\end{equation}

As far as random variables $\zeta_{\e, n}, n = 1, 2, \ldots$ are
concerned, they are interpreted as so-called, ``flag variables'' and are used to
record events $\{\zeta_{\e, n} = 1 \}$ which we interpret as ``rare'' events.

Let us introduce random variables,
\begin{equation}\label{semika}
\xi_\varepsilon = \sum_{n = 1}^{\nu_\varepsilon} \kappa_{\e, n}, \ {\rm where} \ 
\nu_\varepsilon = \min(n \geq 1: \zeta_{\e, n}  =1).
\end{equation}

A random variable $\nu_\varepsilon$ counts the number of transitions
of the imbedded Markov chain $\eta_{\e,n} $ up to the first occurrence of
 ``rare'' event, while a random variable $\xi_\varepsilon$ can be
interpreted as the first-rare-event time
of the first occurrence of  ``rare'' event for the semi-Markov process
$\eta_\e(t)$. 

We also consider the first-rare-event-time  process,
\begin{equation}\label{semikana}
\xi_\varepsilon(t) = \sum_{n = 1}^{[t\nu_\varepsilon]} \kappa_{\e, n}, \ t \geq 0.
\end{equation}

The objective  of this paper is to  describe class $\mathcal F$ of all possible  {\cd} processes $ \xi_0(t), t \geq 0$, which can appear
in the corresponding functional limit theorem given in the form of the asymptotic relation,
$\xi_\varepsilon(t), t \geq 0 \stackrel{\JJ}{\longrightarrow}  \xi_0(t), t \geq 0$ as $\e \to 0$, and to give necessary and sufficient  
conditions for holding of the above  asymptotic relation  with  the specific (by its finite dimensional distributions)  limiting stochastic process 
$\xi_0(t), t \geq 0$ from class $\mathcal F$.


Here and henceforth,  we  use symbol $\stackrel{d}{\longrightarrow}$  to indicate  convergence in 
distribution for random variables (weak convergence of distribution functions) or stochastic processes (weak convergence of finitely dimensional distributions), symbol  $\stackrel{\PP}{\longrightarrow}$ to indicate convergence of random variables in probability,  and symbol  $\stackrel{\JJ}{\longrightarrow}$ to indicate convergence in Skorokhod J-topology  for real-valued {\cd}  stochastic processes defined on 
time interval $[0, \infty)$.  

We refer to books by Gikhman and Skorokhod (1971), Billingsley (1968, 1999) and Silvestrov (2004) for details concerned  the above form of functional convergence. 

The problems formulated above are solved under three general model assumptions.

Let us introduce the probabilities of occurrence of rare event
during one transition step of the semi-Markov process $\eta_\e(t)$,
$$
p_{\e, i} = \mathsf{P}_i \{ \zeta_{\e, 1} = 1 \}, \ i \in \XX.
$$

Here and henceforth, $\mathsf{P}_i$ and $\mathsf{E}_i$ denote,
respectively, conditional probability and expectation calculated
under condition that $\eta_{\e, 0} = i$. 

The first model assumption {\bf A}, imposed on probabilities $p_{i
\varepsilon}$, specifies interpretation of the event $\{\zeta_{\e, n} = 1\}$ as ``rare'' and guarantees the possibility for such
event to occur:
\begin{itemize}
\item[\bf A: ] $0 < \max_{1 \leq i \leq m} p_{\e, i} \to
0$ as $\varepsilon \to 0$.
\end{itemize}

Let us introduce random variables,
\begin{equation}\label{statikani}
\mu_{\e, i}(n) = \sum_{k =1}^n I(\eta_{\e, k-1} = i), \ n = 0, 1, \ldots, \ i \in \XX.
 \end{equation}

If, the Markov chain $\eta_{\e, n}$ is ergodic, i.e., $\mathbb{X}$ is one class of communicative states for  this Markov chain, then  
its stationary distribution is given by the following ergodic relation, 
\begin{equation}\label{statika}
\frac{\mu_{\e, i}(n)}{n}  \stackrel{\PP}{\longrightarrow} \pi_{\e, i}  \ {\rm as} \ n  \to \infty, \  {\rm for} \  i \in \XX.
 \end{equation}
 
The ergodic relation (\ref{statika}) holds for any initial distribution $\bar{q}_{\e}$,  and the stationary distribution 
$\pi_{\e, i}, i \in {\mathbb X}$  does not depend on the initial distribution. Also, all stationary probabilities are positive, i.e., $\pi_i(\e)  > 0, i \in \XX$.

As is known, the stationary probabilities $\pi_i(\e), i \in \XX$ are the unique solution for the system of linear equations,
\begin{equation}\label{statikabas}
\pi_{\e, i} = \sum_{j \in \XX} \pi_{\e, j}  p_{\e, ji}, Êi \in \XX,  \ \ \sum_{i \in \XX} \pi_{\e, i} = 1.
\end{equation}

The second model assumption {is a condition of asymptotically uniform ergodicity for
the embedded Markov chains $\eta_{\e, n}$:
\begin{itemize}
\item[$\mathbf{B}$:] There exists a ring chain of states $i_0, i_1, \ldots, i_N = i_0$ which contains all states from the phase space $\XX$ and such that 
$\varliminf_{\e \to 0}p_{\e, i_{k-1} i_k} > 0$, for $k = 1, \ldots, N$. 
\end{itemize}

As follows from Lemma 1 given below, condition $\mathbf{B}$ guarantees that there exists $\e'_0 \in (0, \e_0]$ such that the Markov chain $\eta_{\e, n}$ is ergodic for every 
$\e \in (0, \e'_0]$. However, condition $\mathbf{B}$ does not require convergence of transition probabilities and, in sequel,  do not imply convergence of stationary probabilities  for the Markov chains $\eta_{\e, n}$ as $\e \to 0$. 

In the case, where the transition probabilities $p_{\e, ij} = p_{0, ij}, i, j \in \XX$ do not depend on parameter $\e$, condition $\bf B$ reduces to the standard  assumption that  the Markov chain $\eta_{0, n}$ with the matrix of transition probabilities $\|  p_{0,ij} \|$ is ergodic.  

Lemma 1 formulated below gives a more detailed information about  condition {\bf B}. 

Finally, the following condition guarantees that the last summand $\kappa_{\e, \nu_\e}$
in the random sum $\xi_\varepsilon$ is asymptotically negligible:
\begin{itemize}
\item[\bf C: ] $\PP_i \{ \kappa_{\e, 1} > \delta / \zeta_{\e, 1} = 1 \} \to 0$ as $\e \to 0$, for $\delta > 0, i \in \XX$.
\end{itemize}

Let us define a probability which is the result of averaging of the
probabilities of occurrence of rare event in one transition step by
the stationary distribution of the imbedded Markov chain $\eta_{\e, n}$,
\begin{equation}\label{defra}
p_\e = \sum\limits_{i=1}^m \pi_{\e, i} p_{\e, i} \ \ {\rm and} \ \ v_\e = p_\e^{-1}.
\end{equation}

Let us  introduce the distribution functions of a sojourn times
$\kappa_{\e,1}$ for the semi-Markov processes $\eta_\e(t)$,
$$
G_{\e, i}(t) = \mathsf{P}_i \{\kappa_{\e, 1} \leq t \}, \ t \geq 0, \ i \in \XX.
$$

Let  $\theta_{\e, n}, n = 1, 2, \ldots$ be i.i.d. random variables
with distribution $G_\e(t)$, which is a result of averaging of
distribution functions of sojourn times by the stationary
distribution of the imbedded Markov chain  $\eta_{\e, n}$,
$$
G_\e(t)=\sum_{i=1}^m \pi_{\e, i} G_{\e, i}(t), \ t \geq 0.
$$

Now, we can  formulate the necessary and sufficient
condition for convergence in distribution for
first-rare-event times: 
\begin{itemize}
\item[$\mathbf{D}$:] $\theta_{\e} = \sum_{n = 1}^{[v_\e]} \theta_{\e, n} \stackrel{d}{\longrightarrow}  \theta_0$ as $\e \to 0$, where $\theta_0$ is a non-negative  random 
variable with  distribution not concentrated in zero.
\end{itemize} 

As well known, {\bf (d$_1$)}   the limiting random variable $\theta_0$ penetrating condition  $\mathbf{D}$ should be infinitely divisible and, thus, its Laplace 
transform has the form,  $\EE e^{-s \theta_0}  = e^{- A(s)}$, where $A(s) = gs + \int_0^\infty (1 - e^{-sv})G(dv), s \geq 0$, $g$ is a non-negative constant and $G(dv)$ is a measure on interval $(0, \infty)$ such that $\int_{(0, \infty)} \frac{v}{1+ v}G(dv) < \infty$;  {\bf (d$_2$)}  $g + \int_{(0, \infty)} \frac{v}{1+ v}G(dv) > 0$ (this is equivalent to the assumption 
that $\PP \{\xi_0 = 0 \} < 1$). 
\vspace{1mm}

Let also consider the homogeneous step-sum process with independent increments (summands are i.i.d. random variables),
\begin{equation}\label{bgre}
\theta_{\e}(t) = \sum_{n = 1}^{[tv_\e]} \theta_{\e, n}, t \geq 0.
\end{equation}

As is known (see, for example, Skorokhod (1964, 1986)), condition $\mathbf{D}$ is necessary and sufficient for holding of the asymptotic relation,
\begin{equation}\label{bgref}
\theta_{\e}(t) = \sum_{n = 1}^{[tv_\e]} \theta_{\e, n}, t \geq 0 \stackrel{\JJ}{\longrightarrow}  \theta_0(t), t \geq 0 \ {\rm as} \ \e \to 0,
\end{equation}
where  $\theta_0(t), t \geq 0$ is a nonnegative L\'{e}vy process (a {\cd} 
homogeneous process with independent increments) with the Laplace transforms $\EE e^{-s \theta_0(t)} = e^{- tA(s)}, s, t  \geq 0$. 

Let us define the Laplace transforms,
$$
\varphi_{\e, i}(s) = \EE_i e^{- s \kappa_{\e, 1}}, i \in \XX, \  \varphi_\e(s) = \EE e^{- s \theta_{\e, 1}} = \sum_{\ \in \XX} \pi_{\e, i} \varphi_{\e, i}(s), \ s \geq 0. 
$$

Condition $\mathbf{D}$ can be reformulated (see, for example, Feller (1966, 1971))  in the equivalent form,  in terms of  the above  Laplace transforms:
\begin{itemize}
\item[$\mathbf{D}_1$:] $v_\e (1 - \varphi_\e(s)) \to A(s)$ as $\e \to 0$, for $s > 0$, where the limiting function $A(s) > 0$, for $s > 0$ and $A(s) \to 0$ as $s \to 0$.
\end{itemize} 

In this case, {\bf (d$_3$)} $A(s)$ is a cumulant of non-negative  random variable with distribution not concentrated in zero. Moreover, {\bf (d$_4$)} $A(s)$ should  be the cumulant of infinitely divisible distribution of the form given in the above conditions {\bf (d$_1$)} and {\bf (d$_2$)}.

The following condition, which is a  variant of the so-called central criterion of convergence (see, for example, Lo\`eve (1977)), is equivalent to condition $\mathbf{D}$,  with the Laplace transform of 
the limiting random variable  $\theta_0$ given in the above conditions {\bf (d$_1$)} and {\bf (d$_2$)}: 
\begin{itemize}
\item[$\mathbf{D}_2$:] {\bf (a)} $v_\e (1-G_\e(u)) \to 
G(u)$ as $\varepsilon \to 0$ for all $u > 0$, which are points
of continuity of the limiting function, which is nonnegative, non-increasing, and
right continuous function defined on interval $(0, \infty)$, with the limiting value $G(+\infty) = 0$;  {\bf (a)} function $G(u)$ is  connected with the measure $G(dv)$ 
 by the relation $G((u', u'']) = G(u') - G(u'')$, $0 < u' \leq u'' <\infty$;
{\bf (b)} $v_\e \int_{(0, u]} v G_\e(dv) \to g + \int_{(0, u]}  v G(dv)$ as $\varepsilon \to 0$ for some
$u > 0$ which is a point of continuity of $G(u)$.
\end{itemize}

It is useful to note that {\bf (d$_5$)}  the asymptotic relation penetrating condition   $\mathbf{D}_2$ {\bf (b)} holds, under condition $\mathbf{D}_2$ {\bf (a)},  for 
any $u > 0$ which is a point of continuity for function $G(u)$.   

In what follows, we also always assume that  asymptotic relations for random variables and processes, defined on trajectories of 
Markov renewal processes $(\eta_{\e, n}, \kappa_{\e, n}, \zeta_{\e, n})$,  hold for any initial distributions $\bar{q}_\e$, if such distributions 
are not specified. 

The main result of the paper is the following theorem. \vspace{1mm}

{\bf Theorem 1.} {\it Let conditions {\bf A}, {\bf B} and {\bf C} hold. Then, 
{\bf (i)} condition {\bf D} is necessary and sufficient for 
holding {\rm (}for some or any initial distributions $\bar{q}_\e$, respectively, in statements of necessity and  sufficiency{\rm )} of the asymptotic relation 
$\xi_\varepsilon = \xi_\varepsilon(1) \stackrel{d}{\longrightarrow}  \xi_0$ as $\e \to 0$, where $\xi_0$ is a non-negative random variable with distribution not concentrated in zero.
In this case, {\bf (ii)} the limiting random variable  $\xi_0$ has the Laplace transform $\EE e^{-s \xi_0} = \frac{1}{1+A(s)}$, where $A(s)$ is a cumulant of infinitely divisible 
distribution defined in  condition  {\bf D}. Moreover, {\bf (iii)} the stochastic processes $\xi_\e(t), t \geq 0 \stackrel{\JJ}{\longrightarrow} \xi_0(t) = \theta_0(t \nu_0), t \geq 0$ as $\e \to 0$, where  {\bf (a)} $\nu_0$ is a random variable, which has  the exponential distribution with parameter $1$, {\bf (b)} $\theta_0(t), t \geq 0$ is a  nonnegative L\'{e}vy  process  with the Laplace transforms $\EE e^{-s \theta_0(t)} = e^{- tA(s)}, s, t  \geq 0$,  {\bf (c)} the random variable $\nu_0$ and 
the process  $\theta_0(t), t \geq 0$  are independent.} \vspace{1mm}

{\bf Remark 1}. According Theorem 1,  class  $\mathcal F$ of all possible nonnegative, nondecreasing,  \cd, stochastically continuous processes $\xi_0(t), t \geq 0$ with distributions of random variables  $\xi_0(t), t > 0$ not concentrated in zero, and such that   the asymptotic relation,  $\xi_\e(t), t \geq 0 \stackrel{\JJ}{\longrightarrow} \xi_0(t), t \geq 0$ as $\e \to 0$,  holds,  coincides  with the class of limiting processes described in proposition {\bf (iii)}. Condition {\bf D} is necessary and sufficient condition for  holding not only the asymptotic relation given in propositions {\bf (i)} -- {\bf (ii)} but also for the much stronger asymptotic relation given in proposition {\bf (iii)}.  \vspace{1mm}

{\bf Remark 2}. The statement ``for some or any initial distributions $\bar{q}_\e$, respectively, in statements of necessity and  sufficiency'' used in the formulation of Theorem1 should be understood in the sense that the asymptotic relation penetrating proposition {\bf (i)} should hold for at least one family of initial distributions $\bar{q}_\e, \e \in (0, \e_0]$, in the  statement of necessity,  and for any family of initial distributions $\bar{q}_\e, \e \in (0, \e_0]$,  in the  statement of sufficiency. \\

{\bf 3. Asymptotics of step-sum reward processes}. \\

Let us consider, for every $\e \in (0, \e_0]$,   the step-sum stochastic process,
\begin{equation}\label{edibas}
\kappa_\e(t) = \sum_{n = 1}^{[tv_\e]} \kappa_{\e, n},  t \geq 0.
\end{equation}

The random variables $\kappa_\e(t) $ can be interpreted as  rewards  accumulated on trajectories of 
the Markov chain  $\eta_{\e, n}$.  Respectively, random variables $\xi_\e$ can be interpreted as rewards  
accumulated on trajectories of the Markov chain  $\eta_{\e, n}$ till the first occurrence of the ``rare''  event.

Asymptotics of the step-sum reward  processes $\kappa_\e(t), t \geq 0$ have its own value. At the same, the corresponding result formulated below 
in Theorem 2 plays the key role in the proof of Theorem 1.

It is useful to note that the flag variables $\zeta_{\e, n}$ are not involved in the  definition of the processes $\kappa_\e(t)$. This let us  replace function $v_\e = p_\e^{-1}$ by an arbitrary function $0 < v_\e \to \infty$ as $\e \to 0$ in condition {\bf D},  Theorem 2  and Lemmas 2 -- 6 formulated below.  
\vspace{1mm}

{\bf Theorem 2.} {\em Let condition {\bf B}  holds. Then, {\bf (i)} condition {\bf D} is necessary and sufficient condition 
for holding {\rm (}for some or any initial distributions $\bar{q}_\e$, respectively, in statements of necessity and  sufficiency{\rm )} of  the asymptotic relation, 
$\kappa_\varepsilon(1) \stackrel{d}{\longrightarrow}  \theta_0$ as $\e \to 0$, where $\theta_0$ is  
a non-negative  random variable with  distribution not concentrated in zero. In this case, {\bf (ii)}  the random variable  $\theta_0$ has the infinitely divisible 
distribution with the Laplace transform $\EE e^{-s \theta_0} = e^{- A(s)}, s \geq 0$ with the cumulant $A(s)$ defined in  
condition  {\bf D}. Moreover, {\bf (iii)} stochastic processes  $\kappa_\e(t), t \geq 0 \stackrel{\JJ}{\longrightarrow} \theta_0(t), t \geq 0$ as $\e \to 0$, where 
$\theta_0(t), t \geq 0$ is a nonnegative  L\'{e}vy  process  with the Laplace transforms 
$\EE e^{-s \theta_0(t)} = e^{- tA(s)}, s, t \geq 0$.}  
\vspace{1mm}

{\bf Remark 3}. According Theorem 2,  class  $\mathcal G$ of all possible nonnegative, nondecreasing,  {\cd},   stochastically continuous 
processes $\theta_0(t), t \geq 0$ with distributions of random variables  $\theta_0(t), t > 0$ not concentrated in zero, and such that the asymptotic relation, $\kappa_\e(t), t \geq 0 \stackrel{\JJ}{\longrightarrow} \theta_0(t), t \geq 0$ as $\e \to 0$,  holds,  coincides 
with the class of limiting processes described in proposition {\bf (iii)}. Condition {\bf D} is necessary and sufficient condition for  holding the asymptotic 
relation given in propositions {\bf (i)} -- {\bf (ii)} as well as for the much stronger asymptotic relation given in proposition {\bf (iii)}.  

We use several useful lemmas in  the proof of Theorems 1 and 2.  \vspace{1mm}

Let $\tilde{\eta}_{\e, n}$ be, for every $\e \in (0, \e_0]$ a Markov chain with the phase space $\XX$ and a matrix of transition 
probabilities $\| \tilde{p}_{\e, ij} \|$. 

We shall use the following condition:
\begin{itemize}
\item[$\mathbf{E}$:] $p_{\e, ij} - \tilde{p}_{\e, ij} \to 0$ as $\e \to 0$, for $i, j \in \XX$.
\end{itemize}

If transition probabilities $\tilde{p}_{\e, ij} \equiv  p_{0, ij}, i,j \in \XX$ do not depend on $\e$, then condition $\mathbf{E}$ reduces to the  following condition:
\begin{itemize}
\item[$\mathbf{F}$:] $p_{\e, ij}\to p_{0, ij}$ as $\e \to 0$, for $i, j \in \XX$.
\end{itemize}

{\bf Lemma 1.} {\em Let  condition $\mathbf{B}$  holds for the Markov chains $\eta_{\e, n}$. Then, 
{\bf (i)} There exists $\e'_0 \in (0, \e_0]$ such that the Markov chain 
$\eta_{\e, n}$ is ergodic, for every $\e \in (0, \e'_0]$ and $0 < \varliminf_{\e \to 0} \pi_{\e, i} \leq  \varlimsup_{\e \to 0} \pi_{\e, i} < 1$, for $i \in \XX$. 
{\bf (ii)} If, together with $\mathbf{B}$, condition  $\mathbf{E}$ holds, then, there exists $\e''_0 \in (0, \e'_0]$ such that Markov chain $\tilde{\eta}_{\e, n}$ is ergodic, for  every  $\e \in (0, \e''_0]$,  and  its stationary distribution $\tilde{\pi}_{\e, i}, i \in \XX$ satisfy the asymptotic relation, $\pi_{\e, i} - \tilde{\pi}_{\e, i} \to 0$ as $\e \to 0$, for $i \in  \XX$.
{\bf (iii)} If condition $\mathbf{F}$  holds, then matrix  $\| p_{0, ij} \|$ is stochastic, condition  $\mathbf{B}$ is equivalent to the assumption that a Markov chain $\eta_{0, n}$, with the matrix of transition probabilities $\| p_{0, ij} \|$,  is ergodic and the following asymptotic relation holds, $\pi_{\e, i} \to \pi_{0, i}$ as $\e \to 0$, for $i \in  \XX$, where $\pi_{0, i},  i \in \XX$ is the stationary distribution of the Markov chain  $\eta_{0, n}$. }
\vspace{1mm}

{\bf Proof}. Let us first prove proposition {\bf (iii)}. Condition $\mathbf{F}$ obviously implies that  matrix $\| p_{0, ij} \|$ is stochastic. Conditions $\mathbf{B}$ and  $\mathbf{F}$ imply that 
$\lim_{\e \to 0}p_{\e, i_{k-1} i_k}  = p_{0, i_{k-1} i_k}  > 0, k = 1, \ldots, N$,  for the ring chain penetrating condition 
$\mathbf{B}$. Thus, the Markov chain $\eta_{0, n}$ with the matrix of transition probabilities $\| p_{0, ij} \|$ is ergodic. Vise versa, the assumption that a Markov chain $\eta_{0, n}$ with the matrix of transition probabilities $\| p_{0, ij} \|$  is ergodic implies that there exists a ring chain of states $i_0, \ldots, i_N = i_0$ which  contains all states from the phase space $\XX$ and such that
$ p_{0, i_{k-1} i_k}  > 0, k = 1, \ldots, N$. In this case, condition $\mathbf{F}$ implies that $\lim_{\e \to 0} p_{\e, i_{k-1} i_k} = p_{0, i_{k-1} i_k}  > 0, k = 1, \ldots, N$, and,  thus, condition {\bf B} holds. Let us assume that the convergence relation for stationary distributions penetrating proposition {\bf (iii)} does 
not hold. In this case, there exist $\delta > 0$ and a sequence $0 < \e_n \to 0$ as $n \to \infty$ such that
$\varliminf_{n \to \infty} |\pi_{\e_n, i'} -  \pi_{0, i'}| \geq \delta$, for some $i' \in  \XX$.  Since, the sequences $\pi_{\e_n, i}, n = 1, 2, \ldots, i \in \XX$ are bounded, there exists
a subsequence $0 < \e_{n_k} \to 0$ as $k \to 0$ such that $\pi_{\e_{n_k}, i} \to  \pi'_{0, i}$ as $k \to \infty$, for $i \in \XX$. This relation, condition 
$\mathbf{F}$ and relation  (\ref{statikabas}) imply that numbers $\pi'_{0, i}, i \in \XX$ satisfy the system of linear equation given in (\ref{statikabas}).  This is impossible, 
since inequality $|\pi'_{0, i'} -  \pi_{0, i'}| \geq \delta$ should hold, while  the stationary distribution $\pi_{0, i}, i \in \XX$ is the unique solution of  system (\ref{statikabas}). 

Let us now prove proposition {\bf (i)}. Condition $\mathbf{B}$ obviously implies that there exist $\e'_0 \in (0, \e_0]$ such that   $p_{\e, i_{k-1} i_k}  > 0, k = 1, \ldots, N$,  for the ring chain penetrating condition $\mathbf{B}$, for $\e \in (0, \e'_0]$. Thus, the Markov chain $\eta_{\e, n}$ is ergodic, for every $\e \in (0, \e'_0]$. Let now assume that $\varliminf_{\e \to 0} \pi_{\e, i'} = 0$, for some $i' \in \XX$. In this case,  there exists a sequence $0 < \e_n \to 0$ as $n \to \infty$ such that
$\pi_{\e_n, i'} \to 0$ as $n \to \infty$.  Since, the sequences $p_{\e_n, ij}, n = 1, 2, \ldots, i, j \in \XX$  are bounded, there exists
a subsequence $0 < \e_{n_k} \to 0$ as $k \to 0$ such that $p_{\e_{n_k}, ij} \to p_{0, ij}$ as $k \to \infty$, for $i, j \in \XX$. By proposition  {\bf (iii)}, the matrix $\| p_{0, ij} \|$ 
is stochastic, the Markov chain   $\eta_{0, n}$ with the matrix of transition probabilities $\| p_{0, ij} \|$ is ergodic and  its stationary distribution $\pi_{0, i}, i \in \XX$ satisfies the asymptotic relation, $\pi_{\e_{n_k}, i} \to  \pi_{0, i}$ as $k \to \infty$, for $i \in  \XX$. This is impossible since equality $\pi_{0, i'} = 0$ should hold, while all stationary probabilities  $\pi_{0, i}, i \in \XX$ 
are positive. Thus, $\varliminf_{\e \to 0} \pi_{\e, i} > 0$, for $i \in \XX$. This implies that, also, $\varlimsup_{\e \to 0} \pi_{\e, i} < 1$, for $i \in \XX$, since $\sum_{i \in \XX}\pi_{\e, i} =1$, 
for $\e \in (0, \e'_0]$.

Finally, let us now prove proposition {\bf (ii)}. Conditions $\mathbf{B}$ and $\mathbf{E}$ obviously imply  that 
$\varliminf_{\e \to 0}\tilde{p}_{\e, i_{k-1} i_k}  = \varliminf_{\e \to 0} p_{\e, i_{k-1} i_k}  > 0, k = 1, \ldots, N$,  for the ring chain penetrating condition 
$\mathbf{B}$. Thus, condition  $\mathbf{B}$ holds also for the Markov chains  $\tilde{\eta}_{\e, n}$  and there exist $\e''_0 \in (0, \e'_0]$ such that 
Markov chain $\tilde{\eta}_{\e, n}$ is ergodic, for  every  $\e \in (0, \e''_0]$. Let assume that the convergence relation for stationary distributions penetrating 
proposition {\bf (ii)} does not hold. In this case,  there exist here exist $\delta > 0$ and a sequence $0 < \e_n \to 0$ as $n \to \infty$ such that
$\varliminf_{n \to \infty} |\pi_{\e_n, i'} -  \tilde{\pi}_{\e_n, i'}| \geq \delta$, for some $i' \in  \XX$. Since, the sequences $p_{\e_n, ij}, n = 1, 2, \ldots, i, j \in \XX$  
are bounded, there exists a subsequence $0 < \e_{n_k} \to 0$ as $k \to 0$ such that $p_{\e_{n_k}, ij} \to p_{0, ij}$ as $k \to \infty$, for $i, j \in \XX$. This relations 
and condition $\mathbf{E}$ imply that, also,   $\tilde{p}_{\e_{n_k}, ij} \to p_{0, ij}$ as $k \to \infty$, for $i, j \in \XX$. By proposition  {\bf (iii)}, the matrix $\| p_{0, ij} \|$ 
is stochastic, the Markov chain   $\eta_{0, n}$ with the matrix of transition probabilities $\| p_{0, ij} \|$ is ergodic and  its stationary distribution $\pi_{0, i}, i \in \XX$ 
satisfies the asymptotic relations, $\pi_{\e_{n_k}, i} \to  \pi_{0, i}$ as $k \to \infty$, for $i \in  \XX$ and $\tilde{\pi}_{\e_{n_k}, i} \to  \pi_{0, i}$ as $k \to  \infty$, for $i \in  \XX$. This is impossible, since  relation $\varliminf_{k \to \infty} |\pi_{\e_{n_k}, i'} -  \tilde{\pi}_{\e_{n_k}, i'}| \geq \delta$ should hold. $\Box$ \vspace{1mm} 

Due to Lemma 1, the asymptotic relation penetrating condition $\mathbf{D}_1$ can, under conditions $\mathbf{A}$, $\mathbf{B}$  and  $\mathbf{E}$,  be  rewritten in the equivalent form, where the stationary probabilities $\pi_{\e, i}, i \in \XX$ are replaced by the stationary probabilities  $\tilde{\pi}_{\e, i}, i \in \XX$, 
\begin{align}\label{bokerwa}
v_\e (1 - \varphi_\e(s)) & = \sum_{i \in \XX} \pi_{\e, i} v_\e(1 -\varphi_{\e, i}(s))   \nonumber \\ 
& \sim  \sum_{i \in \XX} \tilde{\pi}_{\e, i} v_\e(1 -\varphi_{\e, i}(s))    \to A(s) \ {\rm as} \ \e \to 0, \ {\rm for} \ s > 0.
\end{align}

Here and henceforth relation  $a(\varepsilon) \sim b(\varepsilon)$ as $\varepsilon \to 0$ means that $a(\varepsilon)/b(\varepsilon) \to
1$ as $\varepsilon \to 0$.

Proposition {\bf (iii)} of Lemma 1  implies that, in the case, where the transition probabilities $p_{\e, ij} = p_{0,ij}, i, j \in \XX$ do not depend on parameter $\e$ or $p_{\e, ij} \to p_{0, ij}$ as $\e \to 0$, for $i, j \in \XX$,  condition $\bf B$ reduces to the standard  assumption that  the Markov chain $\eta_{0, n}$, with the matrix of transition probabilities 
$\|  p_{0,ij}\|$,  is ergodic.  

These simpler variants of asymptotic ergodicity condition, based on condition $\mathbf{F}$ and the assumption of ergodicity of the Markov chain $\eta_{0, n}$  combined with averaging of characteristic in condition {\bf D} by its stationary distribution $\pi_{0, i},  i \in \XX$,  have been used in the mentioned above works by Silvestrov and Drozdenko (2006a) and Drozdenko (2007a) for proving analogues of propositions {\bf (i)} and {\bf (ii)} of Theorem 1.  

In this case, the averaging of characteristics in the necessary and sufficient condition {\bf D}, in fact, relates mainly to distributions of sojourn times. 
Condition  {\bf B}, used in the present paper, balances in a natural way averaging of characteristics in condition {\bf D} between distributions of sojourn times and stationary distributions of the corresponding embedded Markov chains. 

Let us introduce random variables, which are sequential moments of hitting state $i \in \XX$ by the Markov chain  $\eta_{\e, n}$,
\begin{equation}\label{rewopi}
\tau_{\e, i, n} = \left\{ 
\begin{array}{ll}
\min(k \geq 0, \eta_{\e, k} = i) & \text{for} \  n = 1,  \\
\min(k > \tau_{\e, i, n - 1}, \eta_{\e, k} = i) & {\rm for} \ n \geq 2.
\end{array}
\right.
\end{equation}

Let also define random variables, 
\begin{equation}\label{rewopilo}
\kappa_{\e, i, n} = \kappa_{\e, \tau_{\e, i, n} +1}, n = 1, 2, \ldots, i \in \XX.
\end{equation}

The following simple lemma describe useful properties of the above family of random variables. \vspace{1mm} 

{\bf Lemma 2}. {\em Let condition {\bf B} holds. Then, for every $\e \in (0, \e'_0]$, {\bf (i)} the random variables $\kappa_{\e, i, n}, n = 1, 2, \ldots, i \in \XX$ are independent;
{\bf (ii)} $\PP \{\kappa_{\e, i, n}  \leq t \} = G_{\e, i}(t), t \geq 0$, for $n = 1, 2, \ldots, i \in \XX$; {\bf (iii)} the following representation takes place 
for process $\kappa_\e(t)$,
\begin{equation}\label{edibastabe}
\kappa_\e(t) = \sum_{n = 1}^{[tv_\e]} \kappa_{\e, n} = \sum_{i \in \XX} \sum_{n = 1}^{\mu_{\e, i}([tv_\e])} \kappa_{\e, i, n}, t \geq 0.  
\end{equation}}

It should be noted that the families of random variables $\langle \mu_{\e, i}(n),  n = 0, 1, \ldots, i \in \XX \rangle$ and  $\langle \kappa_{\e, i, n}, n = 1, 2, \ldots, i \in \XX \rangle$ are not independent. 

In what follows, we, for simplicity,  indicate convergence of {\cd}  processes  in uniform U-topology to continuous processes as convergence in J-topology, since, in this case,  convergence J-topology  is equivalent to convergence in uniform U-topology.
\vspace{1mm}

{\bf Lemma 3}. {\em Let condition {\bf B} hold. Then,} 
\begin{equation}\label{edibasta}
\mu^*_{\e, i}(t) =  \frac{\mu_{\e, i}([tv_\e])}{\pi_{\e, i} v_\e}, t \geq 0 \stackrel{\JJ}{\longrightarrow} \mu_{0, i}(t) = t, t \geq 0 \ {\rm as} \ \e \to 0, \ {\rm for} \  i \in \XX. 
\end{equation}

{\bf Proof}. Let $\alpha_{\e, j} = \min(n > 0:  \eta_{\e, n} = j)$ be the moment of first hitting to the state $j \in \XX$ for the Markov chain $\eta_{\e, n}$. Condition {\bf B}  implies that there exist $p \in (0, 1)$ and $\e_p \in (0, \e_0]$ such that  $\prod_{k = 1}^N p_{\e, i_{k-1} i_k} > p$, for $\e \in (0, \e_p]$. The following inequalities are obvious, 
$\PP_i \{ \alpha_{\e, j}  > kN \} \leq (1 - p)^k, k \geq 1,  i, j \in \XX$, for  $\e \in (0, \e_p]$. These inequalities imply that there exists $K_p \in (0, \infty)$ such that  $\max_{i, j \in \XX} \EE_i \alpha_{\e, j}^2  \leq K_p <  \infty, i, j \in \XX$, for $\e \in (0, \e_p]$. Also, as well known, $\EE_i \alpha_{\e, i}  = \pi_{\e, i}^{-1}, i \in \XX$, for $\e \in (0, \e_p]$.

Let $\alpha_{\e, i, n} = \min(k > \alpha_{\e, i, n - 1} : \eta_{\e, k} = i), n = 1, 2, \ldots$ be sequential moments of hitting to state $i \in \XX$ 
for the Markov chain $\eta_{\e, n}$  and $\beta_{\e, i, n} =  \alpha_{\e, i, n} - \alpha_{\e, i, n-1},  n = 1, 2, \ldots$, where  $\alpha_{\e, i, 0} = 0$. The random variables $\beta_{\e, i, n}, n \geq 1$ are independent and  identically 
distributed for  $n \geq 2$. The above relations for moments of random variables $\alpha_{\e, i}  $ imply that $\alpha_{\e, i, 1}/v_\e   \stackrel{\PP}{\longrightarrow} 0$ as $\e \to 0$, for $i \in \XX$.  Also, 
$\PP_i \{ v_\e^{-1} | \alpha_{\e, i, [tv_\e]} - \pi_{\e, i}^{- 1}[tv_\e]|  > \delta \} \leq t K_p/ \delta^2 v_\e, \delta > 0, t \geq 0, i \in \XX$, for $\e \in (0, \e_p]$. These relations obviously implies that random variables 
$\alpha_{\e, i, [tv_\e]}/ \pi_{\e, i}^{-1}v_\e  \stackrel{\PP}{\longrightarrow} t$ as $\e \to 0$, for $t \geq 0$. The dual identities  $\PP \{\mu_{\e, i}(r) \geq k \} = 
\PP \{ \alpha_{\e, i, k} \leq r \}, r, k = 0, 1, \ldots$ let one, in standard way, convert the latter asymptotic relation to the equivalent relation $\mu^*_{\e, i}(t)  = \mu_{\e, i, [tv_\e]}/ \pi_{\e, i}v_\e  \stackrel{\PP}{\longrightarrow} t$ as $\e \to 0$, for $t \geq 0$. Since the processes $\mu^*_{\e, i}(t), t \geq 0$ are nondecreasing and the corresponding limiting function is continuous, the latter asymptotic relation  is (see, for example, Lemma 3.2.2  from Silvestrov (2004))  equivalent to the asymptotic relation (\ref{edibasta}) given in Lemma 3. $\Box$ \vspace{1mm} 

Let now introduce step-sum processes with independent increments,
\begin{equation}\label{edibastaha}
\tilde{\kappa}_\e(t) = \sum_{i \in \XX} \sum_{n = 1}^{[t\pi_{\e, i} v_\e]} \kappa_{\e, i, n}, t \geq 0.  
\end{equation}

Lemmas 2 and 3 let us presume that processes $\tilde{\kappa}_\e(t)$ can be good approximations for processes $\kappa_\e(t)$. \vspace{1mm}

{\bf Lemma 4}. {\em Let condition {\bf B} hold. Then, {\bf (i)} condition {\bf D} holds if and only if the following relation 
holds, $\tilde{\kappa}_\e(1)   \stackrel{d}{\longrightarrow} \theta_0$ as $\e \to 0$, where $\theta_0$ is  
a non-negative  random variable with  distribution not concentrated in zero. In this case, {\bf (ii)}  the random variable  $\theta_0$ has the infinitely divisible 
distribution with the Laplace transform $\EE e^{-s \theta_0} = e^{- A(s)}, s \geq 0$ with the cumulant $A(s)$ defined in  
condition  {\bf D}. Moreover, {\bf (iii)} stochastic processes  $\tilde{\kappa}_\e(t), t \geq 0 \stackrel{\JJ}{\longrightarrow} \theta_0(t), t \geq 0$ as $\e \to 0$, where $\theta_0(t), t \geq 0$ is a nonnegative L\'{e}vy  process with the Laplace transforms 
$\EE e^{-s \theta_0(t)} = e^{- tA(s)}, s, t \geq 0$.} \vspace{1mm}

{\bf Proof of Theorem 2}. The proof of Lemma 4 is an integral  part of the proof of Theorem 2. 

Let us, first, prove that condition {\bf D} implies holding of the asymptotic relations penetrating Lemma 4 and Theorem 2.

Let $\hat{\eta}_{\e, n}, n = 1, 2, \ldots$ be, for every $\e \in (0, \e'_0]$,  a sequence of random variables such that: (a) it is independent of  the Markov chain $(\eta_{\e, n}, \kappa_{\e, n}), n = 0, 1, \ldots$ and (b) it is a sequence of i.i.d. random variables taking value $i$ with probability $\pi_{\e, i}$, for $i \in \XX$.

Note that, in this case,  the sequence of random variables  $\hat{\eta}_{\e, n}, n = 1, 2, \ldots$ is also independent of the  families of random variables $\langle \mu_{\e, i}(n),  n = 0, 1, \ldots, i \in \XX \rangle$ and  $\langle \kappa_{\e, i, n}, n = 1, 2, \ldots, i \in \XX \rangle$.

Let us define random variables, 
\begin{equation}\label{edibasa}
\hat{\mu}_{\e, i}(n) =  \sum_{k = 1}^{n} I(\hat{\eta}_{\e, n} = i),  n = 0, 1, \ldots, i \in \XX.
\end{equation}
and stochastic processes 
\begin{equation}\label{edibahar}
\hat{\kappa}_\e(t) = \sum_{i \in \XX} \sum_{n = 1}^{\hat{\mu}_{\e, i}([t v_\e])} \kappa_{\e, i, n}, t \geq 0.  
\end{equation}
  
Let us also consider the sequence of random variables $\theta_{\e, n} = \kappa_{\e, \hat{\eta}_{\e, n}, n}, n = 1, 2, \ldots$. 
This is the sequence of i.i.d. random variables that follows from the above definition of the sequence of random variables 
$\hat{\eta}_{\e, n}, n = 1, 2, \ldots$ and the family of random variables  
$\kappa_{\e, i, n}, n = 1, 2, \ldots, i \in \XX$. Also,
\begin{equation}\label{edibaheha}
\PP \{\theta_{\e, 1} \leq t \} = \sum_{i \in \XX} \pi_{\e, i} G_{\e, i}(t) = G_\e(t), t \geq 0.
\end{equation}

Let us also define the homogeneous step-sum processes with independent increments using for them, due to relation (\ref{edibaheha}) the 
same notation as for processes  introduced in relation (\ref{bgre}), 
\begin{equation}\label{edibaha}
\theta_\e(t) = \sum_{n = 1}^{[t v_\e]} \theta_{\e, n}, t \geq 0.  
\end{equation}

As well known (see, for example, Skorokhod (1964, 1986)), condition {\bf D} is equivalent to the following relation,
\begin{equation}\label{edibaas}
\theta_\e(t), t \geq 0  \stackrel{d}{\longrightarrow} \theta_0(t), t \geq 0 \ {\rm as} \ \e \to 0. 
\end{equation}

By the definition of the sequence of random variables $\langle \hat{\eta}_{\e, n}, n = 1, 2, \ldots \rangle$ and the family of random variables  
$\langle \kappa_{\e, i, n}, n = 1, 2, \ldots, i \in \XX \rangle$, in particular, due to independence of the above sequence and family,  the following relation holds, 
\begin{equation}\label{edibahawa}
\hat{\kappa}_\e(t), t \geq 0    \stackrel{d}{=} \theta_\e(t), t \geq 0. 
\end{equation}

Relation (\ref{edibahawa}) implies that $\hat{\kappa}_\e(t), t \geq 0$ also is a homogeneous step-sum process with independent increments and 
that condition {\bf D} is equivalent to the following relation,
\begin{equation}\label{edibaastae}
\hat{\kappa}_\e(t), t \geq 0  \stackrel{d}{\longrightarrow} \theta_0(t), t \geq 0 \ {\rm as} \ \e \to 0. 
\end{equation}

Random variables $I(\hat{\eta}_{\e, n} = i), n = 1, 2, \ldots$ are, for every $i \in \XX$,  i.i.d. random variables taking values $1$ and $0$ with probabilities, respectively, 
$\pi_{\e, i}$ and $1 - \pi_{\e, i}$. According proposition {\bf (i)} of Lemma 1, $0 < \varliminf_{\e \to 0} \pi_{\e, i} \leq \varlimsup_{\e \to 0} \pi_{\e, i} < 1$, for every 
$i \in \XX$. Taking into account the  above remarks, this is easy to  prove using the corresponding results from Skorokhod (1964, 1986), that the following relation holds, 
\begin{equation}\label{edibastakamop}
\hat{\mu}^*_{\e, i}(t) = \frac{\hat{\mu}_{\e, i}([tv_\e])}{\pi_{\e, i} v_\e}, t \geq 0 \stackrel{\JJ}{\longrightarrow} \mu_{0, i}(t) = t, t \geq 0 \ {\rm as} \ \e \to 0, \ {\rm for} \  i \in \XX. 
\end{equation}

Let us choose some $0 <  u < 1$.  

By the definition, processes $\tilde{\kappa}_{\varepsilon}(t)$,
$\hat{\kappa}_{\varepsilon}(t)$,  and $\hat{\mu}^*_{\varepsilon, i}(t), i \in \XX$ are
non-negative and non-decreasing. Taking this into account, we get, for $x \geq 0$,
\begin{align}\label{bibermof}
\PP \{\tilde{\kappa}_{\varepsilon}(u) > x \} &  \leq \PP\{ \tilde{\kappa}_{\varepsilon}(u) > x, \hat{\mu}^*_{\e, i}(1)  >  u, i \in \XX \} \nonumber \\
& \quad + \sum_{i \in \XX}\PP\{ \tilde{\kappa}_{\varepsilon}(u) > x, \hat{\mu}^*_{\e, i}(1)  \leq  u \} \nonumber \\
&  \leq \PP \{ \hat{\kappa}_{\varepsilon}(1)  >  x \}  + \sum_{i \in \XX} \PP \{ \hat{\mu}^*_{\e, i}(1)  \leq u \}.
\end{align}

Relations (\ref{edibaastae}), (\ref{edibastakamop}) and inequality (\ref{bibermof}) imply that distributions of random variables $\tilde{\kappa}_{\varepsilon}(u)$ are relatively  compact as $\e \to 0$, 
\begin{align}\label{botyr}
& \lim_{x \to \infty} \varlimsup_{\e \to 0} \PP \{ \tilde{\kappa}_{\varepsilon}(u) > x \}   \leq \lim_{x \to \infty} \varlimsup_{\e \to 0} ( \PP \{\hat{\kappa}_{\varepsilon}(1)  >  x \}  \nonumber \\
& \quad \quad + \sum_{i \in \XX} \PP \{ \hat{\mu}^*_{\e, i}(1)  \leq u \}) =  \lim_{x \to \infty}  \PP \{ \theta_0(1) > x \}  = 0.  
\end{align}

Let also introduce homogeneous step-sum processes with independent increments, for $i \in \XX$, 
\begin{equation}\label{botyra}
\tilde{\kappa}_{\e, i}(t) =  \sum_{n = 1}^{[t\pi_{\e, i} v_\e]} \kappa_{\e, i, n}, t \geq 0.
\end{equation}

Note that, for every $\e \in (0, \e'_0]$, processes $\langle \tilde{\kappa}_{\e, i}(t), t \geq 0 \rangle,  i \in \XX$ are independent.
  
Since, $\tilde{\kappa}_{\e, i}(u) \leq \tilde{\kappa}_{\e}(u)$, for  $i \in \XX$, relation  
(\ref{botyr}) imply that distributions of random variables $\tilde{\kappa}_{\e, i}(1)$ are also relatively  
compact as $\e \to 0$, for every $i \in \XX$, 
\begin{equation}\label{botyras}
\lim_{x \to \infty} \varlimsup_{\e \to 0} \PP \{ \tilde{\kappa}_{\varepsilon, i}(u) > x \}  \leq \lim_{x \to \infty} \varlimsup_{\e \to 0} \PP \{ \tilde{\kappa}_{\varepsilon}(u) > x \} = 0.
\end{equation}
 
This implies that any sequence $0 < \e_n \to 0$ as $n \to \infty$ contains a subsequence $0 < \e_{n_k} \to 0$ as $k \to \infty$ such that random variables,
\begin{equation}\label{botyrasa}
\tilde{\kappa}_{\varepsilon_{n_k}, i}(u) \stackrel{d}{\longrightarrow} \theta_{0, i, u} \ {\rm as} \ k \to \infty, \ {\rm  for} \  i \in \XX,
\end{equation}
where  $\theta_{0, i, u}, i \in \XX$ are  proper nonnegative random variables, with distributions possibly dependent 
of the choice  of subsequence $\e_{n_k}$. 

Moreover, by the central criterion of convergence (see, for example, Lo\`eve (1977)), random variables  $\theta_{0, i, u}, i \in \XX$ have infinitely divisible distributions. Let $\EE e^{- s \theta_{0, i, u}} = e^{- uA_i(s)}, s \geq 0, i \in \XX$ be their Laplace transforms.

As well known (see, for example, Skorokhod (1964, 1986)), relation (\ref{botyrasa})  implies that stochastic processes, 
\begin{equation}\label{botyrame}
\tilde{\kappa}_{\e_{n_k}, i}(t), t \geq 0 \stackrel{\JJ}{\longrightarrow} \theta_{0, i}(t), t \geq 0 \ {\rm as} \ k \to \infty,  \ {\rm  for} \  i \in \XX,
\end{equation}
where  $\theta_{0, i}(t), t \geq 0, i \in \XX$ are nonnegative L\'{e}vy processes with Laplace transforms $\EE e^{- s \theta_{0, i}(t)} = 
e^{- t A_i(s)}, s, t  \geq 0, i \in \XX$, possibly dependent 
of the choice of  subsequence $\e_{n_k}$.  

Moreover, since processes $\tilde{\kappa}_{\e, i}(t), t \geq 0,  i \in \XX$ are independent,  J-conver\-gence of vector 
processes $(\tilde{\kappa}_{\e_{n_k}, 1}(t), \ldots, \tilde{\kappa}_{\e_{n_k}, m}(t)), t \geq 0$ also takes place,
\begin{align}\label{botyrabero}
& (\tilde{\kappa}_{\e_{n_k}, 1}(t), \ldots, \tilde{\kappa}_{\e_{n_k}, m}(t)), t \geq 0 \nonumber \\
& \quad \quad \quad \stackrel{\JJ}{\longrightarrow} (\theta_{0, 1}(t), \ldots, \theta_{0, m}(t)), t \geq 0 \ {\rm as} \ k \to \infty,
\end{align}
where  $\theta_{0, i}(t), t \geq 0, i \in \XX$ are independent nonnegative L\'{e}vy processes with Laplace transforms $\EE e^{- s \theta_{0, i}(t)} = 
e^{- t A_i(s)}, s, t  \geq 0, i \in \XX$, possibly dependent 
of the choice of  subsequence $\e_{n_k}$.  

Note (see, for example, Theorem 3.8.1,  in Silvestrov (2004)) that J-compactness of the vector processes $(\tilde{\kappa}_{\e_{n_k}, 1}(t), \ldots, \tilde{\kappa}_{\e_{n_k}, m}(t))$ follows from J-compactness of their components $\tilde{\kappa}_{\e_{n_k}, i}(t),  i \in \XX,$ since the corresponding limiting processes $\theta_{0, i}(t), i \in \XX$ are stochastically continuous and independent and, thus, they have not with probability $1$ joint points of discontinuity.

Relation (\ref{botyrabero}) obviously implies the following relation,
\begin{equation}\label{botyraberbert} 
\tilde{\kappa}_{\e_{n_k}}(t) = \sum_{i \in \XX} \tilde{\kappa}_{\e_{n_k}, i}(t), t \geq 0 \stackrel{\JJ}{\longrightarrow}  \theta'_{0}(t) = \sum_{i \in \XX}
\theta_{0, i}(t), t \geq 0 \ {\rm as} \ k \to \infty,
\end{equation}
where $\theta_{0, i}(t), t \geq 0, i \in \XX$ are independent nonnegative L\'{e}vy processes described in relation (\ref{botyrabero}).

Since, the limiting processes in (\ref{edibasta}) and (\ref{edibastakamop}) are non-random functions,  relations (\ref{edibasta}), 
(\ref{edibastakamop}) and (\ref{botyraberbert}) imply (see, for example, Subsection 1.2.4 in 
Silvestrov (2004)), by Slutsky theorem, that, 
\begin{align}\label{botyraberik}
& (\mu^*_{\e_{n_k}, 1}(t),  \ldots, \mu^*_{\e_{n_k}, m}(t), \tilde{\kappa}_{\e_{n_k}, 1}(t), \ldots, \tilde{\kappa}_{\e_{n_k}, m }(t)), t \geq 0 \nonumber \\ 
& \quad \quad \stackrel{d}{\longrightarrow} (\mu_{0, 1}(t), \ldots, \mu_{0, m}(t), \theta_{0, 1}(t), \ldots, \theta_{0, m}(t)), t \geq 0 \ {\rm as} \ k \to \infty, 
\end{align}
and
\begin{align}\label{botyraberit}
& (\hat{\mu}^*_{\e_{n_k}, 1}(t), \ldots, \hat{\mu}^*_{\e_{n_k}, m}(t), \tilde{\kappa}_{\e_{n_k}, 1}(t), \ldots, \tilde{\kappa}_{\e_{n_k}, m}(t) ), t \geq 0 \nonumber \\ 
& \quad \quad \stackrel{d}{\longrightarrow} (\mu_{0, 1}(t), \ldots, \mu_{0, m}(t), \theta_{0, 1}(t), \ldots, \theta_{0, m}(t)), t \geq 0 \ {\rm as} \ k \to \infty, 
\end{align}
where $\mu_{0, i}(t) = t, t \geq 0, i \in \XX$ and $\theta_{0, i}(t), t \geq 0, i \in \XX$ are independent nonnegative L\'{e}vy processes defined in 
relation (\ref{botyrabero}).

We can now apply Theorem 3.8.2, from Silvestrov (2004), which give conditions of J-convergence for vector compositions of {\cd} stochastic processes, and get the following asymptotic relations,
\begin{align}\label{botyrabernabaw}
& (\tilde{\kappa}_{\e_{n_k}, 1}(\mu^*_{\e_{n_k}, 1}(t)), \ldots, \tilde{\kappa}_{\e_{n_k}, m}(\mu^*_{\e_{n_k}, m}(t))), t \geq 0  \nonumber \\
& \quad \quad \stackrel{\JJ}{\longrightarrow}  (\theta_{0, 1}(\mu_{0, 1}(t) ), \ldots, \theta_{0, m}(\mu_{0, m}(t) )) \nonumber \\
& \quad \quad \quad =  (\theta_{0, 1}(t), \ldots, \theta_{0, m}(t)), t \geq 0 \ {\rm as} \ k \to \infty, 
\end{align}
and
\begin{align}\label{botyrabernavav}
& (\tilde{\kappa}_{\e_{n_k}, 1}(\hat{\mu}^*_{\e_{n_k}, 1}(t)), \ldots, \tilde{\kappa}_{\e_{n_k}, m}(\hat{\mu}^*_{\e_{n_k}, m}(t))), t \geq 0   \nonumber \\
& \quad \quad \stackrel{\JJ}{\longrightarrow}  (\theta_{0, 1}(\mu_{0, 1}(t) ), \ldots, \theta_{0, m}(\mu_{0, m}(t) )) \nonumber \\
& \quad \quad \quad =  (\theta_{0, 1}(t), \ldots, \theta_{0, m}(t)), t \geq 0 \ {\rm as} \ k \to \infty, 
\end{align}
where $\theta_{0, i}(t), t \geq 0, i \in \XX$ are independent nonnegative L\'{e}vy processes defined in 
relation (\ref{botyrabero}).

Relations (\ref{botyrabernabaw}) and  (\ref{botyrabernavav})  obviously imply J-convergence for sum of components of the 
processes in these relations, i.e. that, respectively, the following relations hold, 
\begin{align}\label{botyrabernaop}
& \kappa_{\e_{n_k}}(t) =  \sum_{i \in \XX} \tilde{\kappa}_{\e_{n_k}, i}(\mu^*_{\e_{n_k}, i}(t)), t \geq 0   \nonumber \\
& \quad \quad \stackrel{\JJ}{\longrightarrow} \theta'_{0}(t) = \sum_{i \in \XX}\theta_{0, i}(t), t \geq 0 \ {\rm as} \ k \to \infty, 
\end{align}
and 
\begin{align}\label{botyrabernaopc}
& \hat{\kappa}_{\e_{n_k}}(t) =  \sum_{i \in \XX} \tilde{\kappa}_{\e_{n_k}, i}(\hat{\mu}^*_{\e_{n_k}, i}(t)), t \geq 0   \nonumber \\
& \quad \quad \stackrel{\JJ}{\longrightarrow} \theta'_{0}(t) = \sum_{i \in \XX}\theta_{0, i}(t), t \geq 0 \ {\rm as} \ k \to \infty, 
\end{align}
where $\theta_{0, i}(t), t \geq 0, i \in \XX$ are independent nonnegative L\'{e}vy processes defined in 
relation (\ref{botyrabero}).

Relation (\ref{edibaastae}) implies that
\begin{equation}\label{botyraberna}
\theta'_{0}(t), t \geq 0 \stackrel{d}{=} \theta_{0}(t), t \geq 0,
\end{equation}

Thus, the limiting process $\theta'_{0}(t) = \sum_{i \in \XX}\theta_{0, i}(t), t \geq 0$ has the  same finite dimensional 
distributions for all  subsequences $\e_{n_k}$ described above. Moreover, the cumulant $A(s)$ of the limiting L\'{e}vy process $\theta_{0}(t)$ is connected with cumulants $A_i(s), i \in \XX$ of L\'{e}vy processes $\theta_{0, i}(t)$ by relation, $A(s) =  \sum_{i \in \XX} A_i(s)$, $s \geq 0$.

Therefore, relations (\ref{botyraberbert}),  (\ref{botyrabernaop}) and  (\ref{botyrabernaopc}) imply that, respectively, the following relations hold,
\begin{equation}\label{botyrabernakolb}
\tilde{\kappa}_{\e}(t) = \sum_{i \in \XX} \tilde{\kappa}_{\e, i}(t), t \geq 0  \stackrel{\JJ}{\longrightarrow} 
\theta_{0}(t), t \geq 0 \ {\rm as} \ \e \to 0,
\end{equation}
and
\begin{align}\label{botyrabernanop}
\kappa_{\e}(t) =  \sum_{i \in \XX} \tilde{\kappa}_{\e, i}(\mu^*_{\e, i}(t)), t \geq 0 
\stackrel{\JJ}{\longrightarrow} \theta_{0}(t), t \geq 0 \ {\rm as} \ \e \to 0, 
\end{align}
as well as, 
\begin{align}\label{botyrabernanok}
\hat{\kappa}_{\e}(t) =  \sum_{i \in \XX} \tilde{\kappa}_{\e, i}(\hat{\mu}^*_{\e, i}(t)), t \geq 0 
\stackrel{\JJ}{\longrightarrow} \theta_{0}(t), t \geq 0 \ {\rm as} \ \e \to 0. 
\end{align}

It is useful to note that relation (\ref{botyrabernanok}) for homogeneous step-sum processes  $\hat{\kappa}_{\e}(t)$ 
follows directly from relation (\ref{edibaastae}). 

It was obtained in the way described above just in order to prove that the limiting process in relations (\ref{botyraberbert}),  
(\ref{botyrabernaop}) and (\ref{botyrabernaopc})   is the same and does not depend on the choice of subsequences $\e_{n_k}$  described above. 
This  made it possible to write down relations (\ref{botyrabernakolb}) and (\ref{botyrabernanop}).

Let us now prove that the asymptotic relation given in proposition {\bf (i)} of Theorem 2 or in proposition {\bf (i)} of Lemma 4 implies
condition {\bf D} to hold.

In both cases, the first step is to prove that distributions of random variables $\tilde{\kappa}_{\varepsilon}(u)$ 
are relatively  compact as $\e \to 0$, for some $u > 0$. 

Let us choose some $0 < u < 1$. 

By the definition, the processes $\kappa_{\varepsilon}(t)$,
$\tilde{\kappa}_{\varepsilon}(t)$,  and $\mu^*_{\varepsilon, i}(t), i \in \XX$ are
nonnegative and nondecreasing.  Taking this into account, we get, for any  $x \geq 0$, 
\begin{align}\label{bibermofm}
\PP \{\tilde{\kappa}_{\varepsilon}(u) > x \} &  \leq \PP\{ \tilde{\kappa}_{\varepsilon}(u) > x, \mu^*_{\e, i}(1)  >  u, 
i \in \XX \} \nonumber \\
& \quad + \sum_{i \in \XX}\PP\{ \tilde{\kappa}_{\varepsilon}(u) > x, \mu^*_{\e, i}(1)  \leq  u \} \nonumber \\
&  \leq \PP \{\kappa_{\varepsilon}(1)  >  x \}  + \sum_{i \in \XX} \PP \{ \mu^*_{\e, i}(1)  \leq u \}.
\end{align}

The asymptotic relation given in proposition {\bf (i)} of Theorem 2, relation (\ref{edibasta}) and inequality (\ref{bibermofm}) imply that, 
\begin{align}\label{botyrvitytbn}
\lim_{x \to \infty} \varlimsup_{\e \to 0} \PP \{ \tilde{\kappa}_{\varepsilon}(u) > x \}  & \leq \lim_{x \to \infty} \varlimsup_{\e \to 0} ( \PP \{\kappa_{\varepsilon}(1)  >  x \}  \nonumber \\
& \quad + \sum_{i \in \XX} \PP \{ \mu^*_{\e, i}(1)  \leq u \}) =  \lim_{x \to \infty}  \PP \{ \theta_0 > x \}  = 0.  
\end{align}

Note that,  in this nessessity  case, the asymptotic relation given in proposition {\bf (i)} of Theorem 2   is 
required  to hold only for at least one family initial distributions $\bar{q}_\e, \e \in (0, \e_0]$.

The asymptotic relation given in proposition {\bf (i)} of Lemma 4  implies that,
\begin{align}\label{botyrvitytb}
\lim_{x \to \infty} \varlimsup_{\e \to 0} \PP \{ \tilde{\kappa}_{\varepsilon}(u) > x \} &  \leq
\lim_{x \to \infty} \varlimsup_{\e \to 0} \PP \{ \tilde{\kappa}_{\varepsilon}(1) > x \} \nonumber \\
&  =  \lim_{x \to \infty}  \PP \{ \theta_0 > x \}  = 0.  
\end{align}

Relation (\ref{botyrvitytbn}), as well as relation (\ref{botyrvitytb}), implies that distributions of random variables 
$\tilde{\kappa}_{\varepsilon}(u)$ are relatively  compact as $\e \to 0$.

Now, we can repeat the part of the above prove related to relations (\ref{botyra}) -- (\ref{botyrabernaopc}). 

Relation (\ref{botyrabernaop}) and the asymptotic relation given in proposition {\bf (i)} of Theorem 2, 
as well as relation (\ref{botyraberbert}) and the asymptotic relation given in proposition {\bf (i)} of Lemma 4,  implies that the random variables $\theta'(1)$ and $\theta_0$, which appears in the above asymptotic relations,
have the same distribution, 
\begin{equation}\label{botyrvitva}
\theta'(1)  \stackrel{d}{=}  \theta_0. 
\end{equation}

Moreover,  cumulant $A(s)$ of the limiting L\'{e}vy process $\theta'_{0}(t)$ coincides with the cumulant of the  random variable 
$\theta_0$, which, therefore, has infinitely divisible distribution. Moreover, relation (\ref{botyrabernaopc}) implies that  cumulant $A(s)$  is connected with cumulants 
$A_i(s), i \in \XX$ of L\'{e}vy processes $\theta'_{0, i}(t)$ by relation $A(s) =  \sum_{i \in \XX} A_i(s), s \geq 0$.

Thus, the limiting process $\theta'_{0}(t), t \geq 0 = \sum_{i \in \XX}\theta_{0, i}(t), t \geq 0$ has the  same finite dimensional 
distributions for all subsequences $\e_{n_k}$ described above.

This let us again to write down relations (\ref{botyrabernakolb}) -- (\ref{botyrabernanok}).

Relation  (\ref{botyrabernanok}) proves, in this case, that condition {\bf D} holds. 

Relation (\ref{botyrabernakolb}) proves proposition {\bf (iii)} of Lemma 4. 

Relation (\ref{botyrabernanop}) proves 
proposition {\bf (iii)} of Theorem 2. $\Box$ \vspace{1mm}

Let us consider the particular case of the model with random variables $\kappa_{\e, n} = f_{\e, \eta_{\e, n-1}}, n = 1, 2, \ldots, i \in \XX$, where $f_{\e, i} \geq 0, i \in \XX$ are nonrandom 
nonnegative numbers. In this case, stochastic process,
\begin{equation}\label{bopser}
\kappa_\e(t) = \sum_{n = 1}^{[tv_\e]} f_{\e, \eta_{\e, n-1}}, t \geq 0.
\end{equation}

Also, the Laplace transforms,
$$
 \varphi_{\e, i}(s) = \EE_i e^{- s f_{\e, \eta_{\e, 0}}} = e^{- s f_{\e, i}}, s \geq 0, \ {\rm for} \ i \in \XX.
$$
and, 
$$
 \varphi_{\e}(s) = \sum_{i \in \XX} \pi_{\e, i} e^{- s f_{\e, i}}, s \geq 0.
$$

Condition $\mathbf{D}_1$ takes, in this case,  the form of the following relation,
\begin{align}\label{njiou}
v_\e (1 - \varphi_\e(s)) & = \sum_{i \in \XX} \pi_{\e, i}v_\e (1 - e^{- s f_{\e}(i)}) \nonumber \\ 
& \to A(s) \ {\rm  as} \  \e \to 0, \ {\rm for} \ s > 0, 
\end{align}
where the limiting function $A(s) > 0$, for $s > 0$ and $A(s) \to 0$ as $s \to 0$. 
 
This condition obviously implies that $1 - \varphi_{\e, i}(s) \to 0$ as $\e \to 0$, for $s >0, i \in \XX$ that is equivalent to relation $f_{\e, i} \to 0$ as $\e \to 0$, 
for $i \in \XX$. In this case, $1 - \varphi_{\e, i}(s) = s f_{\e, i} + o(sf_{\e, i})$ as $\e \to 0$, for every $s > 0, i \in \XX$. These relations let us reformulate 
condition $\mathbf{D}_1$ in terms of functions,
$$
f_\e = v_\e \sum_{i \in \XX} \pi_{\e, i} f_{\e, i}. 
$$

Condition $\mathbf{D}_1$ is equivalent to the following condition:
\begin{itemize}
\item[$\mathbf{G}$:] $f_\e  \to f_0 \in (0, \infty) $ as $\e \to 0$. 
\end{itemize}

Moreover, in this case the cumulant $A(s) = f_0 s, s \geq 0$.

Theorem 2 takes in this case the following form. \vspace{1mm} 

{\bf Lemma  5.} {\em Let condition {\bf B}  holds. Then, {\bf (i)} condition {\bf G} is necessary and sufficient condition 
for holding {\rm (}for some or any initial distributions $\bar{q}_\e$, respectively, in statements of necessity and  sufficiency{\rm )} of  the asymptotic relation, 
$\kappa_\varepsilon(1) \stackrel{d}{\longrightarrow}  \theta_0$ as $\e \to 0$, where $\theta_0$ is  
a non-negative  random variable with  distribution not concentrated in zero. In this case, {\bf (ii)}  the random variable  $\theta_0 \stackrel{d}{=}  f_0$, i.e.,  it is a  constant.  Moreover, {\bf (iii)} stochastic processes  $\kappa_\e(t), t \geq 0 \stackrel{\JJ}{\longrightarrow} f_0 t, t \geq 0$ as $\e \to 0$.}  
\vspace{1mm}

Let us assume that function $f_\e$ satisfy the following natural assumption:
\begin{itemize}
\item[$\mathbf{H}$:] There exists $\e''_0 \in (0, \e'_0]$ such that $f_\e > 0$ for $\e \in (0, \e''_0]$. 
\end{itemize}

In this case, we can describe asymptotic behavior of reward step-sum processes $\kappa_\e(t)$ under weaker than $\mathbf{G}$ condition, which 
admits extremal behavior of functions $f_\e$:
\begin{itemize}
\item[$\mathbf{I}$:] $f_\e  \to f_0 \in [0, \infty] $ as $\e \to 0$. 
\end{itemize}

The following lemma generalizes and supplements Lemma 5.  
\vspace{1mm} 

{\bf Lemma  6.} {\em Let conditions {\bf B} and $\mathbf{H}$  hold. Then, 
 {\bf (i)} $f_\e^{-1} \kappa_\varepsilon(t), t  \geq 0 \stackrel{\JJ}{\longrightarrow}  g_0(t) = t, t \geq 0$ as $\e \to 0$.  {\bf (ii)} Condition {\bf I} is necessary and sufficient condition 
for holding {\rm (}for some or any initial distributions $\bar{q}_\e$, respectively, in statements of necessity and  sufficiency{\rm )} of  the asymptotic relation, 
$\kappa_\varepsilon(1) \stackrel{d}{\longrightarrow}  \theta_0$ as $\e \to 0$, where $\theta_0$ is  
a non-negative   proper or improper random variable. In this case, {\bf (iii)}  the random variable  $\theta_0 \stackrel{d}{\longrightarrow}  f_0$, i.e., it is a constant, and
{\bf (iv)} $\kappa_\varepsilon(t) \stackrel{\PP}{\longrightarrow}  f_0 t$ as $\e \to 0$, for every $t > 0$. Moreover,
{\bf (v)} if $f_0 \in [0, \infty)$ then, $\kappa_\varepsilon(t), t  \geq  0 \stackrel{\JJ}{\longrightarrow}  f_0 t, t \geq 0$ as $\e \to 0${\rm ;}  and {\bf (vi)} if $f_0 = \infty$ then, 
$\min(T, \kappa_\varepsilon(t)), t  > 0 \stackrel{\JJ}{\longrightarrow}  h_T(t) = T, t > 0$ as $\e \to 0$, for every $T > 0$. }  \vspace{1mm}
  
 {\bf Proof}. We can use the following representation,
 \begin{equation}\label{baesd}
 \kappa_\varepsilon(t)  = \sum_{i \in \XX} \mu^*_{\e, i}(t) v_\e \pi_{\e, i} f_{\e, i}, t \geq 0. 
 \end{equation} 

For any sequence $0 < \e_n \to 0$ as $n \to \infty$, there exists a subsequence $0 < \e_{n_k} \to 0$ as $k \to \infty$ such that
\begin{equation}\label{baesda}
\frac{v_{n_k} \pi_{\e_{n_k}, i} f_{\e_{n_k}, i}}{f_{\e_{n_k}}} \to g_i \in [0, 1] \ {\rm as} \ k \to \infty, \ {\rm for} \ i \in \XX.
\end{equation} 

Constants $g_i, i \in \XX$ can depend on the choice of subsequence $\e_{n_k}$, but, obviously satisfy the following relation, 
\begin{equation}\label{baevesda}
\sum_{i \in \XX} g_i = 1.
\end{equation} 

Since the limiting processes in relations  (\ref{edibasta}) given in Lemma 3 are nonrandom functions, relations  (\ref{edibasta}) and  (\ref{baesda}) 
obviously imply that 
\begin{equation}\label{baesdobane}
f_{\e_{n_k}}^{-1} \kappa_{\varepsilon_{n_k}}(t), t \geq 0  \stackrel{d}{\longrightarrow} \sum_{i \in \XX} t g_i =  t, t \geq 0 \ {\rm as} \ k \to \infty.
 \end{equation} 
 
Moreover,  since the processes on the left hand side of the above relation are nondecreasing and the limiting function
is continuous, the following relation (see, for example, Lemma 3.2.2  from Silvestrov (2004)) holds,
\begin{equation}\label{baesdo}
f_{\e_{n_k}}^{-1} \kappa_{\varepsilon_{n_k}}(t), t \geq 0  \stackrel{\JJ}{\longrightarrow} \sum_{i \in \XX} t g_i =  t, t \geq 0 \ {\rm as} \ k \to \infty.
 \end{equation} 
 
Since the limiting process is the same for all subsequences $\e_{n_k}$ described above, relation (\ref{baesdo}) implies that the following relation holds, 
\begin{equation}\label{baesdog}
f_{\e}^{-1} \kappa_{\varepsilon}(t), t \geq 0   \stackrel{\JJ}{\longrightarrow}  g_{0}(t) = t, t \geq 0 \ {\rm as} \ \e \to 0.
\end{equation}

Relation (\ref{baesdog}) implies that random variables $f_{\e}^{-1} \kappa_{\varepsilon}(t)  \stackrel{d}{\longrightarrow} t$ as $\e \to 0$, for every $t \geq 0$.
This implies that the random variables $\kappa_{\varepsilon}(1) = f_\e  \cdot (f_{\e}^{-1} \kappa_{\varepsilon}(1))$ can converge in distribution 
if and only if  $f_\e \to f_0 \in [0, \infty]$ as $\e \to 0$. Moreover, in this case, the limiting (possibly improper) random variable is constant $f_0$, and $\kappa_{\varepsilon}(t) \stackrel{\PP}{\longrightarrow} 
f_0 t$ as $\e \to 0$, for every $t \geq 0$.

If $f_0 \in [0, \infty)$, then the asymptotic relation penetrating propositions {\bf (v)}  can be obtained by application of Theorem 3.2.1 from Silvestrov (2004) 
to processes $\kappa_{\varepsilon}(t) = g_\e(f_{\e}^{-1} \kappa_{\varepsilon}(t)), t \geq 0$, which are 
compositions processes $f_{\e}^{-1} \kappa_{\varepsilon}(t), t \geq 0$ and functions $g_\e(t) = f_\e t, t \geq 0$.

If $f_0  = \infty$ then the asymptotic relation penetrating proposition {\bf (vi)}   can be obtained by application of Theorem 3.2.1 from Silvestrov (2004). to processes $\min(T, \kappa_{\varepsilon}(t)) = h_{\e, T}(f_{\e}^{-1} \kappa_{\varepsilon}(t)), t > 0$, which are  compositions processes $f_{\e}^{-1} \kappa_{\varepsilon}(t), t > 0$ and functions $h_{\e, T}(t) = \min(T, f_\e t), t > 0$. 

Let us now assume that  the asymptotic relation penetrating proposition {\bf (ii)} holds but condition {\bf I} does not hold. 

Relation $f_\e \not\to f_0 \in [0, \infty]$ as $\e \to 0$ holds if and only if 
there exist at least two subsequences $0 < \e'_n, \e''_n \to 0$ as $n \to \infty$ such that (a) $f_{\e'_n} \to f'_0 \in [0, \infty]$ as $n  \to \infty$, (b) $f_{\e''_n} \to f''_0 \in [0, \infty]$ 
as $n  \to \infty$ and (c) $f'_0 \neq f''_0$. In this case, $\kappa_{\e'_n}(1)  \stackrel{\PP}{\longrightarrow} f'_0 $  as $n \to \infty$ and 
$\kappa_{\e''_n}(1)  \stackrel{\PP}{\longrightarrow} f''_0$  as $n \to \infty$ and, thus, random variables $\kappa_{\varepsilon}(1)$ do not converge in 
distribution.    $\Box$ \\

{\bf 4. Asymptotics of first-rare-event times for Markov chains}. \\

The following lemma describe asymptotics for first-rare-event times $\nu_\e$ for Markov chains $\eta_{\e, n}$. 

Note that in this section, we always use function $v_\e = p_\e^{-1}$. 
\vspace{1mm}

{\bf Lemma 7}. {\em Let conditions {\bf A} and {\bf B} hold. Then, the random variables $\nu^*_\e = p_\e \nu_\e \stackrel{d}{\longrightarrow} \nu_0$ as $\e \to 0$, where 
$\nu_0$ is a random variable exponentially distributed with parameter $1$.} 
\vspace{1mm}

{\bf Proof}. Let us define probabilities, for $\e \in (0, \e_0]$,  
$$
P_{\e, i j} = \PP_i \{ \eta_{\e, 1} = j, \zeta_{\e, 1} = 0 \}, \
 \tilde{p}_{\e, i j} = \frac{P_{\e, i j}}{\sum_{r \in \XX} P_{\e, i r}} = \frac{P_{\e, i j}}{1 - p_{\e, i}}, \ i, j \in \XX.
$$

Let also $\tilde{\eta}_{\e, n}, n = 0, 1, \ldots$ be a homogeneous Markov chain with the phase space $\XX$, an initial distribution $\bar{q}_\e = \langle q_{\e, i}, i \in \XX \rangle$ and the matrix
of transition probabilities $\|  \tilde{p}_{\e, i j}  \|$. 

The following relation  takes place, for  $t \geq 0$,
\begin{equation*}
\begin{aligned}\label{burte}
\PP \{ \nu^*_\e > t \} & = \sum_{i \in \XX} q_{\e, i} \sum_{i = i_0, i_1, \ldots, i_{[tv_\e]} \in \XX} \prod_{k = 1}^{[tv_\e]} P_{\e, i_{k-1} i_k} \nonumber \\
& = \sum_{i \in \XX} q_{\e, i} \sum_{i = i_0, i_1, \ldots, i_{[tv_\e]} \in \XX} \prod_{k = 1}^{[tv_\e]} \tilde{p}_{\e, i_{k-1} i_k} (1 - p_{\e, i_{k-1}}) \nonumber \\
\end{aligned}
\end{equation*}
\begin{align}\label{burtek}
& = \EE \exp\{ - \sum_{k = 1}^{[tv_\e]} - \ln (1 - p_{\e, \tilde{\eta}_{\e, k -1}}) \}. 
\end{align}

Conditions {\bf A} and {\bf B} imply that condition {\bf B} holds for transition probabilities of the Markov chains $\tilde{\eta}_{\e, n}$, since, 
the following relation holds, for $i, j \in \XX$,
\begin{align}\label{burteb}
| p_{\e, ij} -  \tilde{p}_{\e, i j} | & = \frac{|p_{\e, ij}(1 - p_{\e, i}) - P_{\e, i j}|}{1 - p_{\e, i}} \nonumber \\
& = \frac{|Ê\PP_i \{ \eta_{\e, 1} = j, \zeta_{\e, 1} = 0 \} - p_{\e, ij}p_{\e, i}|}{1 - p_{\e, i}} \nonumber \\
& \leq  \frac{2p_{\e, i}}{1 - p_{\e, i}}  \to 0 \ {\rm as} \ \e \to 0.
\end{align}

Thus, by Lemma 1, there exist $\e''_0 \in (0, \e'_0]$ such that the Markov chain $\tilde{\eta}_{\e, n}$ is ergodic, for every $\e_0 \in (0, \e''_0]$, and its stationary probabilities $\tilde{\pi}_{\e, i}, i \in \XX$ satisfy the following relation,
\begin{equation}\label{nopi}
\tilde{\pi}_{\e, i} - \pi_{\e, i} \to 0 \ {\rm as} \ \e \to 0,  \ {\rm for}  \ i \in \XX.
\end{equation}

We can apply Lemma 5, which is a particular case of Theorem 2,  to the nonnegative step-sum process, 
\begin{equation}\label{olpi}
\kappa^*_{\e}(t) = \sum_{n = 1}^{[tv_\e]} - \ln(1 - p_{\e, \tilde{\eta}_{\e, n-1}}), t \geq 0.
\end{equation}

To do this, we should check that condition {\bf G}  holds for functions  $f_\e(i) = - \ln(1 - p_{\e, i}), i \in \XX$.
Indeed, using condition {\bf A}, {\bf B}, Lemma 1 and relation (\ref{nopi}), we get,
\begin{align}\label{llnu}
f_\e & =  - v_\e \sum_{i \in \XX} \tilde{\pi}_{\e, i} \ln (1 - p_{\e, i}) \sim  v_\e   \sum_{i \in \XX} \tilde{\pi}_{\e, i}  \, p_{\e, i} \nonumber \\ 
& \sim   v_\e \sum_{i \in \XX} \pi_{\e, i} \, p_{\e, i} =  v_\e p_\e = 1 \ {\rm as} \ \e \to 0.
\end{align}

This relation is a variant of condition {\bf G}. In this case the corresponding limiting constant $\theta_0 = 1$  and 
the process $\theta_0(t) = t, t \geq 0$ is a non-random linear function. By applying sufficiency proposition of Lemma 5 to the step-sum 
process $\kappa^*_{\e}(t)$, we get the following relation,
\begin{equation}\label{firsaba}
\kappa^*_{\e}(t), t \geq 0 \stackrel{d}{\longrightarrow} \theta_0(t) = t, t \geq 0  \ {\rm as} \ \e \to 0.
\end{equation}

The expression on the right hand side of relation (\ref{burte}) is, just, the Laplace transform of the nonnegative random variable $\kappa^*_{\e}(t)$ at point $1$. Thus, 
relation (\ref{firsaba}) implies, by continuity theorem for Laplace transforms, that the following relation holds, for every $t \geq 0$,
\begin{equation}\label{burtelo}
\PP \{ \nu^*_\e > t \}  = \EE e^{- \kappa^*_{\e}(t)} \to e^{- t}  \ {\rm as} \ \e \to 0.
\end{equation}

The proof is complete. $\Box$

Let, as in Lemma 8,  $f_{\e, i}, i \in \XX$ be nonrandom nonnegative numbers and $f_\e = v_\e \sum_{i \in \XX} \pi_{\e, i} f_{\e, i}$. Let us introduce stochastic processes, 
\begin{equation}\label{olpil}
\nu_{\e}(t) = \sum_{n = 1}^{[t\nu_\e]} f_{\e, \eta_{\e, n-1}}, t \geq 0.
\end{equation}

The following lemma generalizes Lemma 7 and is used in what follows. \vspace{1mm}

{\bf Lemma  8.} {\em Let conditions {\bf A}, {\bf B} and {\bf H}   hold. Then, 
 {\bf (i)} $f_\e^{-1} \nu_\varepsilon(t), t  \geq 0 \stackrel{\JJ}{\longrightarrow}  t \nu_0, t \geq 0$ as $\e \to 0$, where $\nu_0$ is a random variable exponentially distributed with parameter $1$.  {\bf (ii)} Condition {\bf I} is necessary and sufficient condition 
for holding {\rm (}for some or any initial distributions $\bar{q}_\e$, respectively, in statements of necessity and  sufficiency{\rm )} of  the asymptotic relation, $\nu_\varepsilon(1) \stackrel{d}{\longrightarrow}  \nu$ as $\e \to 0$, where $\nu$ is  
a non-negative  random variable with  distribution not concentrated in zero. In this case, {\bf (iii)}  the random variable  $\nu \stackrel{d}{=} f_0 \nu_0$. Moreover, {\bf (iv)}  if $f_0 \in [0, \infty)$ then, $\nu_\varepsilon(t), t  \geq  0 \stackrel{\JJ}{\longrightarrow}  f_0 \nu_0 t, t \geq 0$ as $\e \to 0$,  and, {\bf (v)} if $f_0 = \infty$ then, $\min(T, \nu_\varepsilon(t)), t  > 0 \stackrel{\JJ}{\longrightarrow}  h_T(t) = T, t > 0$ as $\e \to 0$, for every $T > 0$ and, thus, {\bf (vi)} 
$\nu_\varepsilon(t) \stackrel{\PP}{\longrightarrow}  \infty$ as $\e \to 0$, for  $t > 0$.}  \vspace{1mm}

{\bf Proof}. The following representation takes place,
\begin{equation}\label{olpilhako}
\nu_{\e}(t)  = \kappa_\e(t\nu^*_\e), t \geq 0,
\end{equation}
where $\kappa_\e(t)$ are processes defined in relation (\ref{bopser}).

Relations given in proposition {\bf (i)} of Lemma 6 and in Lemma 7 imply, by Slutsky theorem, the following relation,
\begin{equation}\label{olpilhamr}
(t \nu^*_\e,  f_\e^{-1}  \kappa_\e(t)), t \geq 0  \stackrel{d}{\longrightarrow} (t \nu_0, t), t \geq 0 \ {\rm as } \ \e \to 0.
\end{equation}

The components of the processes on the left hand side of relation (\ref{olpilhamr}) are non-decreasing processes and the process on the right hand side of relation 
(\ref{olpilhamr})  is continuous.  This let us apply Theorem  3.2.1 from Silvestrov (2004) to processes $f_\e^{-1}\nu_{\e}(t) = f_\e^{-1}\kappa_\e(t\nu^*_\e), t \geq 0$ and to get the asymptotic relation penetrating the proposition {\bf (i)} of Lemma  8. 

Relation penetration proposition {\bf (i)} of Lemma 8 implies that random variables $f_\e^{-1} \nu_\varepsilon(1) \stackrel{d}{\longrightarrow} \nu_0$ as $\e \to 0$. This implies  that random variables $\nu_{\varepsilon}(1) = f_\e  \cdot (f_{\e}^{-1} \nu_{\varepsilon}(1))$ can converge in distribution 
if and only if  $f_\e \to f_0 \in [0, \infty]$ as $\e \to 0$. Moreover, in this case, the limiting (possibly improper) random variable $\nu \stackrel{d}{=} f_0 \nu_0$.

If $f_0 \in [0, \infty)$, then relations given in proposition {\bf (iv)} of Lemma 6 and in Lemma 7 imply, by Slutsky theorem, the following relation,
\begin{equation}\label{olpilhamrnok}
(t \nu^*_\e,   \kappa_\e(t)), t \geq 0  \stackrel{d}{\longrightarrow} (t \nu_0,  f_0 t), t \geq 0 \ {\rm as } \ \e \to 0.
\end{equation}

The components of the processes on the left hand side of relation (\ref{olpilhamrnok}) are non-decreasing processes and the process on the right hand side of relation (\ref{olpilhamr})  is continuous.  This let us apply Theorem 3.2.1 from Silvestrov (2004) to processes $\nu_{\e}(t) = \kappa_\e(t\nu^*_\e), t \geq 0$ and to get the asymptotic relation penetrating the proposition {\bf (iv)} of Lemma  8. 

If $f_0 = \infty$, then relations given in proposition {\bf (v)} of Lemma 6 and in Lemma 7 imply, by Slutsky theorem, the following relation,
\begin{equation}\label{olpilhamrno}
(t \nu^*_\e,   \min(T, \kappa_\e(t))), t > 0  \stackrel{d}{\longrightarrow} (t \nu_0,  T), t > 0 \ {\rm as } \ \e \to 0, \ {\rm for} \ T > 0.
\end{equation}

The components of the processes on the left hand side of relation (\ref{olpilhamrno}) are non-decreasing processes and the process on the right hand side of relation (\ref{olpilhamr})  is continuous.  Also the limiting random variable $t \nu_0 > 0$ with probability $1$, for every $t > 0$.  This let us apply Theorem 3.2.1 (and the remarks made in Subsection 3.2.6) from Silvestrov (2004)   to processes $ \min(T, \nu_{\e}(t)) = \min(T, \kappa_\e(t\nu^*_\e)), t > 0$ and to get the asymptotic relation penetrating the proposition {\bf (v)} of Lemma  8. 

Proposition {\bf (vi)} of this lemma is the direct corollary of proposition {\bf (v)}.

Let us now assume that  the asymptotic relation penetrating proposition {\bf (ii)} holds but condition {\bf I} does not hold. 

Relation $f_\e \not\to f_0 \in [0, \infty]$ as $\e \to 0$ holds if and only if 
there exist at least two subsequences $0 < \e'_n, \e''_n \to 0$ as $n \to \infty$ such that (a) $f_{\e'_n} \to f'_0 \in [0, \infty]$ as $n  \to \infty$, (b) $f_{\e''_n} \to f''_0 \in [0, \infty]$  as $n  \to \infty$ and (c) $f'_0 \neq f''_0$. In this case, $\nu_{\e'_n}(1)  \stackrel{d}{\longrightarrow} f'_0 \nu_0 $  as $n \to \infty$ and 
$\nu_{\e''_n}(1)  \stackrel{d}{\longrightarrow} f''_0\nu_0$  as $n \to \infty$ and, thus, random variables $\nu_{\varepsilon}(1)$ do not converge in 
distribution. $\Box$ \\

{\bf 5. Asymptotics of first-rare-event times for semi-Markov \\ \makebox[10mm]{} processes}. \\

{\bf Proof of Theorem 1}. Now we are prepared to complete the proof of this theorem. 
Let us, first, concentrate attention on propositions {\bf (i)} and  {\bf (ii)} of this theorem.

Let us introduce  Laplace transforms, 
$$
\varphi_{\e, i j}( \imath, s) = \EE_i I(\eta_{\e, 1} = j, \zeta_{\e, 1} = \imath ) e^{- s \kappa_{\e, 1}}, s \geq 0, \ {\rm for} \ i, j \in \XX, \ \imath = 0, 1, 
$$
and 
$$
\varphi_{\e, i}(\imath, s) = \EE_i I(\zeta_{\e, 1} = \imath ) e^{- s \kappa_{\e, 1}},  s \geq 0,  \ {\rm for} \ i \in \XX, \ \imath = 0, 1. 
$$

Let also introduce conditional Laplace transforms,
$$
\phi_{\e, i j}( \imath, s) = \EE_i \{ I(\eta_{\e, 1} = j) e^{- s \kappa_{\e, 1}} / \zeta_{\e, 1} = \imath \}, s \geq 0, \ {\rm for} \ i, j \in \XX, \ \imath = 0, 1, 
$$
and
$$
\phi_{\e, i}(\imath, s) = \EE_i  \{ e^{- s \kappa_{\e, 1}} / \zeta_{\e, 1} = \imath \}, s \geq 0,  \ {\rm for} \ i \in \XX, \ \imath = 0, 1.
$$

Now, let  us  define probabilities,  for $s \geq 0$,
$$
p_{\e, s,  i j}  = \frac{\varphi_{\e, i j}(0, s)}{\sum_{r \in \XX} \varphi_{\e, i j}(0, s) } = \frac{\varphi_{\e, i j}(0, s)}{\varphi_{\e, i}(0, s)}, \ i, j \in \XX.
$$

Let  $(\eta_{\e, s, n}, \zeta_{\e, s, n}), n = 0, 1, \ldots$ be, for every $s \geq 0$,  a Markov renewal process, with the phase space $\XX \times  \{0, 1 \}$, the initial an initial distribution $\bar{q}_\e = \langle q_{\e, i} = \PP \{\eta_{\e, 0} = i, \zeta_{\e, s, 0} = 0 \} = \PP \{\eta_{\e, s, 0} = i \}, i \in {\mathbb X} \rangle$  and transition probabilities,
\begin{equation}\label{sadop}
\begin{aligned}
& \PP \{ \eta_{\e, s, n+1} = j, \zeta_{\e, s, n+1} = \jmath
/ \eta_{\e, s, n} = i, \zeta_{\e, s, n} = \imath \}   \\
& \quad \quad = \PP \{ \eta_{\e, s, n+1} = j, \zeta_{\e, s, n+1} = \jmath / \eta_{\e, s, n} = i \}  \\
& \quad \quad = p_{\e, s,  i j}(p_{\e, i} \,  \jmath + (1 - p_{\e, i})(1 - \jmath))  , \ i, j \in \mathbb{X}, \ \imath, \jmath = 0, 1.
\end{aligned}
\end{equation}

Note that the firat component of the Markov renewal process, $\eta_{\e, s, n}, n = 0, 1, \ldots$ is a homogeneous Markov chain with the phase space $\XX$, an initial distribution $\bar{q}_\e = \langle q_{\e, i}, i \in \XX \rangle$ and the matrix of transition probabilities $\| p_{\e, s, i j}  \|$. 

Let us also  introduce random variables, 
\begin{equation}\label{sadopna}
\nu_{\e, s} = \min(n \geq 1: \zeta_{\e, s, n} = 1).
\end{equation}

Let us prove that condition {\bf D} or conditions  {\bf A}, {\bf B}  and the asymptotic relation penetrating proposition {\bf (i)} of Theorem 1 imply 
that, for every $s \geq  0$, condition {\bf B} holds for transition probabilities of the Markov chain $\eta_{\e, s, n}$.  

Condition {\bf D} obviously, implies that, for  $i \in \XX$, 
\begin{equation}\label{gorta}
\varphi_{\e, i}(s) \to 1 \ {\rm  as} \ \e \to 0, \ {\rm for} \  s \geq 0, 
\end{equation}

Let us show that  conditions {\bf A},  {\bf B}  and the asymptotic relation penetrating proposition {\bf (i)} of Theorem 1 also implies that relation (\ref{gorta}) holds. 

Let us use representation,
\begin{equation}\label{gortana}
\xi_\e =  \sum_{i \in \XX} \sum_{n = 1}^{\mu_{\e, i}(\nu_\e)} \kappa_{\e, i, n} =   \sum_{i \in \XX} \sum_{n = 1}^{[\mu^*_{\e, i}(\nu_\e)\pi_{\e, i} v_\e]} \kappa_{\e, i, n}.
\end{equation}

Let us now assume that relation (\ref{gorta}) does not holds. This means that there exists $i \in \XX$ such that for some $\delta, p  > 0$ and $\e_{\delta, p} \in (0, \e'_0]$ 
probability $\PP \{  \kappa_{\e, i, 1} \geq \delta \} \geq p$, for $\e \in (0, \e_{\delta, p}]$. This obviously implies that random variables  $\tilde{\kappa}_{\e, i}(t) = \sum_{n = 1}^{[t \pi_{\e, i}v_\e]} \kappa_{\e, i, n}  \stackrel{\PP}{\longrightarrow} \infty$ as $\e \to 0$,  for $t  > 0$, and, thus, stochastic processes  $\min (T, \tilde{\kappa}_{\e, i}(t)), t > 0
 \stackrel{d}{\longrightarrow} h_T(t) = T, t > 0$ as $\e \to 0$. Since, the processes $\tilde{\kappa}_{\e, i}(t), t > 0$ are non-decreasing and the limiting function 
$h_T(t) = T, t > 0$  is continuous, the latter relation implies (see, for example, Theorem 3.2.1 from Silvestrov (2004)) that  
 $\min (T, \tilde{\kappa}_{\e, i}(t)), t > 0 \stackrel{\JJ}{\longrightarrow} h_T(t) = T, t \geq 0$ as $\e \to 0$.
 Also, by Lemma 8, applied to the model with functions $f_{\e, j} = I( j = i) (\pi_{\e, i} v_\e)^{-1}, j \in \XX$, the following relation takes place, $\mu^*_{\e, i}(\nu^*_\e) \stackrel{d}{\longrightarrow} \nu_0$ as $\e \to 0$, where $\nu_0$ is a random variable exponentially distributed with parameter $1$. The latter two relations imply, by Slutsky theorem, 
that $(\mu^*_{\e, i}(\nu^*_\e), \min (T, \tilde{\kappa}_{\e, i}(t))), t > 0 \stackrel{d}{\longrightarrow} (\nu_0, h_T(t)), t > 0$ as $\e \to 0$. Now we can apply 
Theorem 2.2.1 from Silvestrov (2004) that yields the following relation,   $\min (T, \tilde{\kappa}_{\e, i}(\mu^*_{\e, i}(\nu^*_\e))) \stackrel{d}{\longrightarrow} T$ as 
$\e \to 0$, for any $T > 0$.  This is possible  only if $\tilde{\kappa}_{\e, i}(\mu^*_{\e, i}(\nu^*_\e))) \stackrel{\PP}{\longrightarrow} \infty$ as $\e \to 0$. Thus,  
random variables $\tilde{\kappa}_{\e, i}(\mu^*_{\e, i}(\nu^*_\e))) = \sum_{n = 1}{\mu_{\e, i}(\nu_\e)} \kappa_{\e, i, n}  \leq \xi_\e \stackrel{\PP}{\longrightarrow} \infty$ as $\e \to 0$. This relation contradicts to  the asymptotic relation penetrating proposition {\bf (i)} of \mbox{Theorem 1.} 

Relation (\ref{gorta}) and condition {\bf A} imply the following relation,
\begin{equation}\label{edibasbasno}
\varphi_{\e, i}(s, 0) = \EE_i I(\zeta_{\e, 1} = 0 ) e^{- s \kappa_{\e, 1}}
\to 1 \ {\rm as} \ \e \to 0, \ {\rm for} \ s \geq 0, i, j \in \XX.
\end{equation}
which implies that, for $s \geq 0$,
\begin{align}\label{edibastame}
p_{\e, ij} - p_{\e, s, i j}  & = \frac{p_{\e, ij}\varphi_{\e, i}(0, s)  - \varphi_{\e, i j}(0, s)}{\varphi_{\e, i}(0, s)} \nonumber \\
& \leq \frac{| p_{\e, ij}\varphi_{\e, i}(0, s) - p_{\e, ij}| + |p_{\e, ij} - \varphi_{\e, i j}(0, s)|}{\varphi_{\e, i}(0, s)} \nonumber \\
& \leq \frac{p_{\e, ij}|  \varphi_{\e, i}(0, s) - 1| +\EE_i I(\eta_{\e, 1} = j) |1 -  I(\zeta_{\e, 1} = 0) e^{- \kappa_{\e, 1}})|}{\varphi_{\e, i}(0, s)} \nonumber \\
& \leq \frac{2(1 - \varphi_{\e, i}(0, s))}{\varphi_{\e, i}(0, s)} \to 0 \ {\rm as} \ \e \to 0, \ {\rm for} \  i, j \in \XX.
\end{align}

Thus, for every $s \geq 0$, there exist $\e'_{0, s} \in (0, \e_0]$ such that the Markov chain $\tilde{\eta}_{\e, n, s}$ is ergodic, for every $\e \in (0, \e'_{0, s}]$, and 
its stationary probabilities $\pi_{\e, s, i}, i \in \XX$ satisfy the following relation, 
\begin{equation}\label{nopibol}
\pi_{\e, s, i} - \pi_{\e, i} \to 0 \ {\rm as} \ \e \to 0,  \ {\rm for}  \ i \in \XX.
\end{equation}

Let us assume that Markov chains $\eta_{\e, n}$ and $\eta_{\e, n, s}$ has the same  initial distribution $\bar{q}_\e$.

The following representation  takes place for the Laplace transform of the random variables $\xi_\e$, for $s \geq 0$,    
\begin{align}\label{voptrewopiba}
\EE e^{- s \xi_\e} & = \sum_{i \in \XX} q_{\e, i} \sum_{n = 1}^{\infty}  \sum_{i = i_0, i_1, \ldots, i_{n} \in \XX} 
\prod_{k = 1}^{n-1} \varphi_{\e, i_{k-1} i_k}(0, s) \varphi_{\e, i_{n-1} i_n}(1, s) \nonumber \\
& =\sum_{i \in \XX} q_{\e, i} \sum_{n = 1}^{\infty}  \sum_{i = i_0, i_1, \ldots, i_{n-1} \in \XX} 
\prod_{k = 1}^{n-1} \varphi_{\e, i_{k-1} i_k}(0, s)  \sum_{i_n \in \XX} \varphi_{\e, i_{n-1} i_n}(1, s) \nonumber \\
& = \sum_{i \in \XX} q_{\e, i} \sum_{n = 1}^{\infty}  \sum_{i = i_0, i_1, \ldots, i_{n-1} \in \XX} 
\prod_{k = 1}^{n-1} p_{\e, s, i_{k-1} i_k} \nonumber \\ 
& \quad \times (1 - p_{\e, i_{k-1}}) \phi_{\e, i_{k-1}}(0, s)  p_{\e, i_{n-1}} \sum_{i_n \in \XX} \frac{\varphi_{\e, i_{n-1}, i_n}(1, s)}{p_{\e, i_{n-1}}} \nonumber \\
& = \sum_{i \in \XX} q_{\e, i} \sum_{n = 1}^{\infty}  \sum_{i = i_0, i_1, \ldots, i_{n-1} \in \XX} 
\prod_{k = 1}^{n-1} p_{\e, s, i_{k-1} i_k} \nonumber \\ 
& \quad \times (1 - p_{\e, i_{k-1}}) \phi_{\e, i_{k-1}}(0, s) p_{\e, i_{n-1}} \phi_{\e, i_{n-1}}(1, s) \nonumber \\
& = \EE \exp\{ - \sum_{k = 1}^{\nu_{\e, s}} - \ln \phi_{\e, \eta_{\e, s, k-1}}(0, s) \nonumber \\ 
& \quad - \ln \phi_{\e, \eta_{\e, s, \nu_{\e, s} -1}}(0, s)  +  \ln \phi_{\e, \eta_{\e, s,  \nu_{\e, s} - 1}}(1, s) \}.  
\end{align} 

Relation (\ref{edibasbasno}) and condition {\bf A} imply that the following relation holds, 
\begin{equation}\label{edibasbasbi}
\phi_{\e, i}(0, s) =  \frac{\phi_{\e, i}(0, s)}{1 - p_{\e, i}}
\to 1 \ {\rm as} \ \e \to 0,  \ {\rm for} \ s \geq 0, i \in \XX. 
\end{equation}

Also condition {\bf C} is equivalent to the following relation,  
\begin{align}\label{edibbi}
\phi_{\e, i}(1, s) & =  \EE_i \{  e^{- s \kappa_{\e, 1}} / \zeta_{\e, 1} = 1 \}  \to 1 \ {\rm as} \ \e \to 0, \ {\rm for} \ s \geq 0, i \in \XX. 
\end{align}

The above two relations obviously imply that, 
\begin{equation}\label{edibbifer}
| \ln \phi_{\e, \eta_{\e, s, \nu_{\e, s} -1}}(0, s) | +  | \ln \phi_{\e, \eta_{\e, s,  \nu_{\e, s} - 1}}(1, s)  | \stackrel{\PP}{\longrightarrow} 0 \ {\rm  as} \ \e \to 0, \ {\rm for} \ s \geq  0.  
\end{equation}

Representation (\ref{voptrewopiba}) and relation (\ref{edibbifer}) imply the following relation, 
\begin{equation}\label{edibbiferno}
\EE e^{- s \xi_\e}  \sim  \EE e^{ - \tilde{\nu}_{\e, s}} \ {\rm  as} \ \e \to 0, \ {\rm for} \ s > 0.
\end{equation}
where
\begin{equation}\label{olpivil}
\tilde{\nu}_{\e, s} = \sum_{n = 1}^{\nu_{\e, s}} -  \ln \phi_{\e, \eta_{\e, s, k-1}}(0, s). 
\end{equation}

Relations (\ref{nopibol}), (\ref{edibasbasbi}) and  proposition {\bf (i)}  of Lemma 1  imply that, 
\begin{align}\label{edibbimase}
A_{\e}(s) & = - v_\e \sum_{i \in \XX} \pi_{\e, s, i} \ln \phi_{\e, i}(0, s) \nonumber \\ 
& \sim  v_\e  \sum_{i \in \XX} \pi_{\e, s, i} (1- \phi_{\e, i}(0, s))   \nonumber \\
&\sim v_\e \sum_{i \in \XX} \pi_{\e, i}  (1 - \phi_{\e, i}(0, s))  \ {\rm  as}  \ \e \to 0, \ {\rm for} \ s >  0.
\end{align}

Let us assume that  condition {\bf D}  holds additionally to conditions conditions {\bf A} -- {\bf C}. 

Condition {\bf D} is equivalent to condition {\bf D}$_1$, and, thus, due to relations (\ref{edibasbasbi}) and (\ref{edibbi}), condition {\bf A} and 
proposintion {\bf (i)} of Lemma 1, to the following relation, 
\begin{align}\label{edibbimaseboi}
v_\e(1 - \varphi_\e(s)) & =  v_\e \sum_{i \in \XX} \pi_{\e, i}  (1 - \varphi_{\e, i}(s))   \nonumber \\
& = v_\e \sum_{i \in \XX} \pi_{\e, i} (1 - (1- p_{\e, i}) \phi_{\e, i}(0, s) -  p_{\e, i}\phi_{\e, i}(1, s))  \nonumber \\
& = v_\e \sum_{i \in \XX} \pi_{\e, i} ( (1- p_{\e, i}) (1 - \phi_{\e, i}(0, s)) +  p_{\e, i}(1 - \phi_{\e, i}(1, s))  \nonumber \\
& \sim v_\e \sum_{i \in \XX} \pi_{\e, i}  (1- p_{\e, i}) (1 - \phi_{\e, i}(0, s))   \nonumber \\
& \sim v_\e \sum_{i \in \XX} \pi_{\e, i} (1 - \phi_{\e, i}(0, s)) 
\to A(s) \ {\rm  as}  \ \e \to 0, \ {\rm for} \ s > 0, 
\end{align}
where $A(s) > 0$, for $s > 0$ and $A(s) \to 0$ as $s \to 0$.

Relations (\ref{edibbimase}) and (\ref{edibbimaseboi}) imply that, in this case,   
\begin{align}\label{edibbioi}
A_{\e}(s)  = - v_\e \sum_{i \in \XX} \pi_{\e, s, i} \ln \phi_{\e, i}(0, s)  \to A(s) \ {\rm  as}  \ \e \to 0, \ {\rm for} \ s > 0.
\end{align}

Now, we can, for every $s > 0$, apply the sufficiency statement  of proposition {\bf (iv)} of Lemma 8 to random variables $\tilde{\nu}_{\e, s}$. 
This yields, the following relation, 
\begin{equation}\label{oivil}
\tilde{\nu}_{\e, s}  \stackrel{d}{\longrightarrow}  A(s)\nu_0  \ {\rm  as}  \ \e \to 0, \ {\rm for} \ s > 0, 
\end{equation}
where $\nu_0$ is exponentially distributed random variable with parameter $1$.

This relation implies, by continuity theorem for Laplace transforms, the following relation, 
\begin{equation}\label{edibbifernobomi}
\EE e^{- s \xi_\e}  \sim  \EE  e^{-  \tilde{\nu}_{\e, s}} \to   \EE e^{- A(s)\nu_0} = \frac{1}{1 + A(s)}  \ {\rm  as} \ \e \to 0,  \ {\rm for} \ s > 0.
\end{equation}

Relation (\ref{edibbifernobomi}) proves sufficiency statements of propositions {\bf (i)} and  {\bf (ii)} of Theorem 1.

Let now assume that conditions  {\bf A} -- {\bf C}  plus proposition {\bf (i)} of Theorem 1 hold.

The asymptotic relation (in proposition {\bf (i)} of Theorem 1)  expressed in terms of Laplace transforms takes the form 
of relation (which should be assumed to hold for some initial distributions $\bar{q}_\e$),
\begin{equation}\label{ediifbo}
\EE e^{- s \xi_\e} \to e^{ - A_0(s)} \ {\rm as} \  \e \to 0, \ {\rm for} \ s > 0, 
\end{equation}
where $A_0(s) > 0$ for $s > 0$ and $A_0(s) \to 0$ as $s \to 0$. 

Let us assume that conditions  {\bf A} -- {\bf C}  hold but condition {\bf D}  does not holds. 

The latter assumption  means, due to relation (\ref{edibbimase}), that either  (a) $A_{\e}(s) \to A(s)  \in (0, \infty)$ as $s \to 0$, for every $s>0$, but  
$A(s)  \not\to 0$ as $\e \to 0$, or (b) $A_{\e}(s^*) \not\to A(s^*) \in (0, \infty)$ as $\e \to 0$, for some $s^* > 0$. The latter relation holds if and only if 
there exist at least two subsequences $0 < \e'_n, \e''_n \to 0$ as $n \to \infty$ such that (b$_1$) $A_{\e'_n}(s^*) \to A'(s^*) \in [0, \infty]$ as $n  \to \infty$, (b$_2$)
$A_{\e''_n}(s^*) \to A''(s^*) \in [0, \infty]$ as $n  \to \infty$ and (b$_3$) $A''(s^*) < A'(s^*)$.

In the case (a), we can repeat the part of the above proof presented in relations (\ref{edibbimase}) -- (\ref{edibbifernobo}) and, taking into account relation 
 (\ref{ediifbo}),  to get relation, $\EE e^{- s \xi_\e}  \sim  \EE e^{ - \tilde{\nu}_{\e, s}}  \to \frac{1}{1 + A(s)}  = e^{ - A_0(s)}$  as $\e \to 0$,  for  $s > 0$.
This relation implies that $A(s) \to 0$ as $\e \to 0$, i.e., the case (a) is impossible.

In the case (b), sub-case, $A'(s^*) = \infty$,  is impossible. Indeed, as was shown in the proof of 
Lemma 8, applied to random variables $\tilde{\nu}_{\e, s^*}$,  in this case,  $\tilde{\nu}_{\e'_n, s^*} \stackrel{\PP}{\longrightarrow} \infty$ as $n \to \infty$, and, thus,
$\EE e^{- s^* \xi_{\e'_n}}  \sim  \EE  e^{-  \tilde{\nu}_{\e'_n, s^*}}  \to  0$   as $n \to \infty$. This relation contradicts to relation (\ref{ediifbo}).

Sub-case,  $A''(s^*) = 0$, is also impossible. Indeed, as was shown in the proof of 
Lemma 8, random variables $\tilde{\nu}_{\e, s^*}$, in this case, $\tilde{\nu}_{\e''_n, s^*} \stackrel{\PP}{\longrightarrow} 0$ as $n \to \infty$, and, thus,
$\EE e^{- s^* \xi_{\e''_n}}  \sim  \EE  e^{-  \tilde{\nu}_{\e''_n, s^*}}  \to  1$   as $n \to \infty$. This relation also contradicts to relation (\ref{ediifbo}).

Finally, the remaining sub-case,   $0 < A''(s^*) < A'(s^*) < \infty$,  is also impossible. Indeed,  the sufficiency statement of Lemma 7 applied to random variables 
$\tilde{\nu}_{\e, s^*}$ yields, in this case, two relations  $\tilde{\nu}_{\e'_n, s^*} \stackrel{d}{\longrightarrow} A'(s^*) \nu_0$ as $n \to \infty$ and 
$\tilde{\nu}_{\e''_n, s^*} \stackrel{d}{\longrightarrow} A''(s^*) \nu_0$ as $n \to \infty$, where $\nu_0$ is exponentially distributed random variable with parameter $1$. These 
relations imply that $\EE e^{- s^* \xi_{\e'_n}}  \sim  \EE  e^{-  \tilde{\nu}_{\e'_n, s^*}}  \to  \frac{1}{1 +A'(s^*)}$   as $n \to \infty$ and $\EE e^{- s^* \xi_{\e''_n}}  \sim  
\EE  e^{-  \tilde{\nu}_{\e''_n, s^*}}  \to  \frac{1}{1 +A''(s^*)}$  as $n \to \infty$. These relations contradict to relation (\ref{ediifbo}), since  $\frac{1}{1 +A'(s^*)}  \neq  \frac{1}{1 +A''(s^*)}$.  

Therefore, condition {\bf D}  should hold. This complete the proof of propositions  {\bf (i)} and  {\bf (ii)} of Theorem 1. 

The following lemma brings together the asymptotic relations given in Theorem 2 and Lemma 7. The proposition of this lemma 
gives the last intermediate result required for completing  the proof of proposition  {\bf (iii)} in   Theorem 1. \vspace{1mm}

{\bf Lemma 9.} {\em Let conditions {\bf A}, {\bf B}, {\bf C} and  {\bf D}   hold. Then, the following asymptotic relation holds,
$(\nu^*_\e, \kappa_\varepsilon(t)), t \geq 0 \stackrel{d}{\longrightarrow}  (\nu_0,  \theta_0(t)), t \geq 0$ as $\e \to 0$, {\bf (a)} $\nu_0$ is a random variable, which has  the exponential distribution with parameter $1$, {\bf (b)} $\theta_0(t), t \geq 0$ is a  nonnegative L\'{e}vy  process  with the Laplace transforms $\EE e^{-s \theta_0(t)} = e^{- tA(s)}, s, t  \geq 0$, with the cumulant $A(s)$ defined in  condition  {\bf D}, {\bf (c)} the random variable $\nu_0$ and the process  $\theta_0(t), t \geq 0$  are independent.} 
\vspace{1mm}

{\bf Proof}.  The following representation takes place, for $s, t \geq 0$, 
\begin{align}\label{burtebo}
\EE I( \nu^*_\e > t) e^{- s \kappa_\e(t)} & = \sum_{i \in \XX} q_{\e, i} 
\sum_{i = i_0, i_1, \ldots, i_{[tv_\e]} \in \XX} \prod_{k = 1}^{[tv_\e]} \varphi_{\e, i_{k-1} i_k}(0, s) \nonumber \\
& = \sum_{i \in \XX} q_{\e, i} \sum_{i = i_0, i_1, \ldots, i_{[tv_\e]} \in \XX} \prod_{k = 1}^{[tv_\e]} p_{\e, s,  i_{k-1} i_k} \nonumber \\
& \quad \times (1 - p_{\e, i_{k-1}}) \phi_{\e, i_{k-1}}(0, s)  \nonumber \\
& = \EE \exp\{ - \sum_{k = 1}^{[tv_\e]} (- \ln (1 - p_{\e, \tilde{\eta}_{\e, k -1}})  \nonumber \\
& \quad - \ln \phi_{\e, \tilde{\eta}_{\e, k -1}}(0, s) ) \}. 
\end{align}

Using condition {\bf A}, {\bf B}, Lemma 1 and relation (\ref{nopibol}), we get, for  $s \geq 0$,  the following analogue of relation (\ref{llnu}),
\begin{align}\label{llnumop}
f_{\e, s} & =  - v_\e \sum_{i \in \XX} \tilde{\pi}_{\e, s, i} \ln (1 - p_{\e, i}) \sim  v_\e   \sum_{i \in \XX} \tilde{\pi}_{\e, s, i}  \, p_{\e, i} \nonumber \\ 
& \sim   v_\e \sum_{i \in \XX} \pi_{\e, i} \, p_{\e, i} =  v_\e p_\e = 1 \ {\rm as} \ \e \to 0.
\end{align}

Relations (\ref{edibbioi}) and (\ref{llnumop}) and  imply that Lemma 5 can, for every $s >  0$, be applied to the processes,
\begin{equation}\label{edibbifernobo}
\kappa_{\e, s}(t) = \sum_{k = 1}^{[tv_\e]} (- \ln (1 - p_{\e, \tilde{\eta}_{\e, k -1}}) - \ln \phi_{\e, \tilde{\eta}_{\e, k -1}}(0, s) ), t \geq 0.
\end{equation}

This yields that the following relation holds, for every $s > 0$,
\begin{equation}\label{edibbobov}
\kappa_{\e, s}(t), t \geq 0 \stackrel{d}{\longrightarrow} t + A(s)t, t \geq 0 \ {\rm as} \ \e \to 0.
\end{equation}

Let us denote, for $i, j \in \XX, n = 0, 1, \ldots, s \geq 0$, 
$$
\Psi_{\e, ij}(n, s) = \EE_i I( \nu_\e > n, \eta_{\e, n} = j) e^{- s \sum_{k = 1}^{n} \kappa_{\e, k} },  
$$
and
$$
\Psi_{\e, i}(n, s) = \EE_i I( \nu_\e > n) e^{- s \sum_{k = 1}^{n} \kappa_{\e, k} } = \sum_{j \in \XX} \Psi_{\e, ij}(n, s).  
$$ 

Relation (\ref{edibbobov}) implies, by continuity theorem for Laplace transforms, the following relation, for $t \geq 0$,
\begin{align}\label{burtebonops}
\EE I( \nu^*_\e > t) e^{- s \kappa_\e(t)} & = \Psi_{\e, i}([t v_\e], s)   \nonumber \\ 
& = \EE e^{- \kappa_{\e, s}(t)} \to e^{-t -A(s)t}  \nonumber \\
& =  e^{-t} e^{- A(s)t} \ {\rm as} \ \e \to 0, \ {\rm for} \ s > 0.
\end{align}

Let us also denote, for $i, j \in \XX, n = 0, 1, \ldots, s \geq 0$, 
$$
\psi_{\e, ij}(n, s) = \EE_i I(\eta_{\e, n} = j) e^{-s \sum_{k = 1}^n \kappa_{\e, k}},
$$
and 
$$
\psi_{\e, i}(n, s) =  \EE_i  e^{-s \sum_{k = 1}^n \kappa_{\e, k}} = \sum_{j \in \XX} \psi_{\e, ij}(n, s).
$$

Relation (\ref{burtebonops}) easily  implies  that, for $s > 0$ and $0 \leq t'' \leq t' < \infty$, 
\begin{align}\label{easy}
\Psi_{\e, i}([t' v_\e] - [t'' v_\e], s)  & \sim \Psi_{\e, i}([(t' - t'') v_\e], s) \nonumber \\
& \to  e^{- (t' - t'')} e^{- A(s)(t' - t'')} \ {\rm  as} \  \e \to 0, 
\end{align}

 Also the proposition {\bf (iii)} of Theorem 2 easily implies that, for $s > 0$ and $0 \leq t'' \leq t' < \infty$.
\begin{align}\label{easyd}
\psi_{\e, i}([t' v_\e] - [t'' v_\e], s) &  \sim \psi_{\e, i}([(t' - t'') v_\e], s) \nonumber \\
& \to  e^{- A(s)(t' - t'')} \ {\rm  as} \  \e \to 0.  
\end{align}

Relations (\ref{easy}) and (\ref{easyd}) imply that, for $s >0$ and $0 \leq t'' \leq t' < \infty$,  
\begin{align}\label{edibbobovt}
& \sum_{j \in \XX} \Psi_{\e, i j}([t' v_\e] -  [t'' v_\e], s) \nonumber \\ 
& \quad \quad \quad = \Psi_{\e, i}([t' v_\e] -  [t'' v_\e], s) \to e^{- (t' - t'')} e^{- A(s)(t' - t'')}  \ {\rm  as} \  \e \to 0.
\end{align}
and 
\begin{align}\label{edibbobovger}
& \sum_{j \in \XX} \psi_{\e, i j}([t' v_\e] -  [t'' v_\e], s)  \nonumber \\ 
& \quad \quad \quad = \psi_{\e, i}([t' v_\e] -  [t'' v_\e], s) \to e^{- A(s)(t' - t'')}  \ {\rm  as} \  \e \to 0.
\end{align}

Now, we shall use the following representation for multivariate joint distributions   of random variable $\nu^*_\e$ and increments 
of stochastic process $\kappa_\varepsilon(t)$ for $0 = t_0 \leq  t_1 < \cdots t_k = t \leq t_{k +1} \leq \cdots \leq t_n < \infty, 1 \leq k < n < \infty$ and $s_1, \ldots, s_n \geq 0$,
\begin{equation*}
\begin{aligned} 
& \EE I( \nu^*_\e > t_k)\exp \{ - \sum_{r= 1}^n s_r (\kappa_\e(t_r) - \kappa_\e(t_{r-1}) \} \nonumber \\
& \quad \quad \quad = \sum_{i_0, \ldots, i_{n}  \in \XX} q_{\e, i_0} \prod_{r = 1}^k \Psi_{\e, i_{r-1} i_r}([t_r v_\e] -  [t_{r-1} v_\e], s_r) \nonumber \\
& \quad \quad \quad  \quad  \times \prod_{r = k+1}^n \psi_{\e, i_{r-1} i_r}([t_r v_\e] -  [t_{r-1} v_\e], s_r) \nonumber \\
& \quad \quad \quad = \sum_{i_0  \in \XX} q_{\e, i_0} \sum_{i_1 \in \XX} \Psi_{\e, i_{0} i_1}([t_1 v_\e] -  [t_{0} v_\e], s_1)   \nonumber \\
\end{aligned}
\end{equation*}
\begin{align}\label{burtnopsva}
& \quad \quad \quad  \quad \cdots  \times \sum_{i_{k} \in \XX} \Psi_{\e, i_{k-1} i_k}([t_k v_\e] -  [t_{k-1} v_\e], s_k)  \nonumber \\
& \quad \quad \quad  \quad  \times \sum_{i_{k +1} \in \XX} \psi_{\e, i_{k} i_{k+1}}([t_{k+1} v_\e] -  [t_{k} v_\e], s_{k+1}) \nonumber \\
& \quad \quad \quad  \quad \cdots  \times \sum_{i_{n} \in \XX} \psi_{\e, i_{n-1} i_n}([t_{n} v_\e] -  [t_{n-1} v_\e], s_{n}) 
\end{align}

Using relations (\ref{edibbobovt}),  (\ref{edibbobovger}) and representation (\ref{burtnopsva}) we get recurrently, for $0 = t_0 \leq  t_1 < \cdots t_k = t \leq t_{k +1} \leq \cdots \leq t_n < \infty, 1 \leq k < n < \infty$ and $s_1, \ldots, s_n > 0$,
\begin{align}\label{burtnopges}
& \EE I( \nu^*_\e > t_k)\exp \{ - \sum_{r= 1}^n s_r (\kappa_\e(t_r) - \kappa_\e(t_{r-1})) \} \nonumber \\
& \quad \quad  \sim \EE I( \nu^*_\e > t_k)\exp \{ - \sum_{r= 1}^{n-1} s_r (\kappa_\e(t_r)  - \kappa_\e(t_{r-1})) \} e^{- A(s_n)(t_{n} - t_{n-1})} \nonumber \\
&  \quad \quad \cdots \sim \EE I( \nu^*_\e > t_k)\exp \{ - \sum_{r= 1}^{k} s_r (\kappa_\e(t_r)  - \kappa_\e(t_{r-1})) \} \nonumber \\
&  \quad \quad \quad \quad \times \exp\{ \sum_{r = k+1}^n - A(s_r)(t_{r} - t_{r-1}) \} \nonumber \\
& \quad \quad  \sim \EE I( \nu^*_\e > t_{k-1})\exp \{ - \sum_{r= 1}^{k-1} s_r (\kappa_\e(t_r)  - \kappa_\e(t_{r-1})) \} \nonumber \\
&  \quad \quad \quad \quad \times \exp\{- (t_k - t_{k-1}) \} \exp \{ \sum_{r = k}^n - A(s_r)(t_{r} - t_{r-1}) \} \nonumber \\
& \quad \quad \cdots \sim \exp\{- \sum_{r = 1}^k (t_r - t_{r-1})\} \exp \{ \sum_{r = 1}^n - A(s_r)(t_{r} - t_{r-1}) \} \nonumber \\ 
& \quad \quad = \exp\{- t \} \exp \{ \sum_{r = 1}^n - A(s_r)(t_{r} - t_{r-1}) \} \ {\rm  as} \  \e \to 0.
\end{align}

This relation is equivalent an  form of the asymptotic relation given in Lemma 9. $\Box$ \vspace{1mm}

Now, we can complete the proof of Theorem 1.

The asymptotic relation given in Lemma 9 can, obviously, be rewritten in the following equivalent form,   
\begin{equation}\label{final}
(t\nu^*_\e, \kappa_\varepsilon(t)), t \geq 0 \stackrel{d}{\longrightarrow}  (t\nu_0,  \theta_0(t)), t \geq 0  \ {\rm  as} \  \e \to 0, 
\end{equation}
where the random variable $\nu_0$ and the stochastic process $\theta_0(t), t \geq 0$ are  described in Lemma 9.

Asymptotic relation given in proposition {\bf (iii)} of Theorem 2 and relation (\ref{final}) let us apply Theorem 3.4.1 from Silvestrov (2004) to the compositions of stochastic processes  $\kappa_\varepsilon(t), t \geq 0$ and $t\nu^*_\e, t \geq 0$  that yield the following relation,
\begin{equation}\label{finala}
\xi_\e(t) = \kappa_\varepsilon(t\nu^*_\e), t \geq 0 \stackrel{\JJ}{\longrightarrow}  \theta_0(t\nu_0), t \geq 0  \ {\rm  as} \  \e \to 0.
\end{equation} 

The proof of Theorem 1 is complete. $\Box$

\end{document}